\documentclass[a4paper,12pt]{amsart}
\usepackage{tikz-cd} 
\usepackage{tcolorbox}
\usepackage[justification=centering,labelformat=empty]{caption}
\usepackage{amsopn,amsfonts,amssymb,amsthm,mathtools,amsmath}
\usepackage{color}
\usepackage{dsfont}
\usepackage{latexsym}

\usepackage{wrapfig}
\usepackage[T1]{fontenc}
\usepackage{lmodern}
\usepackage{textcomp}
\usepackage{soul} 
\usepackage{rotating}
\usepackage{xcolor}
\usepackage{wasysym}
\newcommand\doublecheck{\textcolor{black}{\checked\kern-0.6em\checked}}

\usepackage{hyperref}
\hypersetup{colorlinks=true,linkcolor=blue,citecolor=red}
\usepackage{multicol}
\usepackage[utf8]{inputenc}
\usepackage{float}
\usepackage{hhline}
\usepackage{multirow}
\usepackage{anysize}
\marginsize{2cm}{2cm}{2cm}{2cm}
\usepackage{graphicx}

\usepackage{lscape}

\def\g{\mathfrak{g}}
\def\h{\mathfrak{h}}
\def\n{\mathfrak{n}}

\def\s{\mathfrak{s}}

\def\B{\mathcal{B}}

\def\l{\ell}

\def\C{\mathbb{C}}
\def\R{\mathbb{R}}

\def\Z{\mathbb{Z}}
\def\N{\mathbb{N}}

\newcommand{\SL}{\operatorname{SL}}

\def\ad{\operatorname{ad}}

\def\Tr{\operatorname{Tr}}
\def\Id{\operatorname{Id}}

\def\I{\operatorname{I}}

\newcommand{\Aut}{\operatorname{Aut}}

\newcommand{\Der}{\operatorname{Der}}

\def\alt{\raise1pt\hbox{$\bigwedge$}}

\def\la{\langle}
\def\ra{\rangle}
\def\multiset#1#2{\ensuremath{\left(\kern-.3em\left(\genfrac{}{}{0pt}{}{#1}{#2}\right)\kern-.3em\right)}}
\def\nil{\mathfrak{nil}}

\theoremstyle{plain}
\newtheorem{theorem}{\bf Theorem}[section]

\newtheorem{proposition}[theorem]{\bf Proposition}
\newtheorem{lemma}[theorem]{\bf Lemma}

\theoremstyle{definition}

\theoremstyle{remark}
\newtheorem{remark}[theorem]{\bf Remark}

\newcommand{\matriz}[1]{\ensuremath{\begin{bmatrix}#1\end{bmatrix}}}
\newcommand{\smatriz}[1]{\ensuremath{\left[\begin{smallmatrix}#1\end{smallmatrix}\right]}}   

\newcommand{\comillas}[1]{``#1''}

\title[]{Six-dimensional complex solvmanifolds with non-invariant trivializing sections of their canonical bundle} 

\author{Alejandro Tolcachier}

\address{Dipartimento di Scienza ed Alta Tecnologia, Università degli Studi dell'Insubria, Via Valleggio 11,
	22100, Como, Italy}
\email{alejandro.tolcachier@uninsubria.it}


\begin{document}
\begin{abstract}
It is known that there exist complex solvmanifolds $(\Gamma\backslash G,J)$ whose canonical bundle is trivialized by a holomorphic section which is not invariant under the action of $G$. The main goal of this article is to classify the six-dimensional Lie algebras corresponding to such complex solvmanifolds, thus extending the previous work of Fino, Otal and Ugarte for the invariant case. To achieve this, we complete the classification of six-dimensional solvable strongly unimodular Lie algebras admitting complex structures and identify among them, the ones admitting complex structures with Chern-Ricci flat metrics. 
Finally we construct complex solvmanifolds with non-invariant holomorphic sections of their canonical bundle. In particular, we present an example of one such solvmanifold that is not biholomorphic to a complex solvmanifold with an invariant holomorphic section of its canonical bundle. Additionally, we discover a new $6$-dimensional solvable strongly unimodular Lie algebra equipped with a complex structure that has a non-zero holomorphic $(3,0)$-form.

\end{abstract}

\maketitle
\section{Introduction}\label{S: intro}
The canonical bundle of a complex manifold $(M,J)$ of complex dimension $n$, denoted $K_{(M,J)}$, is defined as the $n$-th exterior power of its holomorphic cotangent bundle, and it is a holomorphic line bundle over $M$. This line bundle is holomorphically trivial when there exists a nowhere vanishing $(n,0)$-form which is holomorphic (or equivalently, closed). Complex manifolds with holomorphically trivial canonical bundle, often equipped with a special Hermitian metric, play a relevant role both in geometry and in theoretical physics. For instance, compact K\"ahler manifolds $M^{2n}$ with global Riemannian holonomy contained in $\operatorname{SU}(n)$ have holomorphically trivial canonical bundle. More generally, any Calabi-Yau manifold (i.e., a compact K\"ahler manifold $M$ with $c_1(M)=0$ in $H^2(M,\R)$) has holomorphically torsion canonical bundle, that is, $K_{(M,J)}^{\otimes k}$ is trivial for some $k\in \N$. According to \cite{To}, (non-Kähler) compact complex manifolds with holomorphically torsion canonical bundle have vanishing first Bott-Chern class, $c_1^{BC}=0$, and therefore they are examples of \textit{non-Kähler Calabi-Yau} manifolds.  Special attention is paid to real dimension six, due to theoretical physical applications. Indeed, the solutions of the Hull-Strominger system occur in compact complex manifolds $M$ endowed with a special Hermitian metric (not necessarily K\"ahler) and $K_{(M,J)}$ holomorphically trivial (see for instance \cite{OU}).  An important source of these distinguished manifolds is provided by compact quotients $\Gamma\backslash G$, where $G$ is a simply connected Lie group and $\Gamma$ is a cocompact and discrete subgroup of $G$ (called a \textit{uniform lattice}), equipped with an \textit{invariant} complex structure. For instance, when  $\Gamma\backslash G$ is a nilmanifold (i.e. $G$ is nilpotent), it was shown in \cite{BDV} that the simply connected nilpotent Lie group $G$ admits a non-vanishing left-invariant holomorphic $(n,0)$-form $\sigma$ (with $\dim_\R G=2n$), by using a distinguished basis of left-invariant $(1,0)$-forms provided by Salamon in \cite{Sal}. Since $\sigma$ is left-invariant, it induces an invariant trivializing section of $K_{(\Gamma\backslash G,J)}$ for any lattice $\Gamma \subset G$. 

In contrast, it is known that for complex solvmanifolds (i.e. $G$ is solvable) several different phenomena can occur. There are examples of complex solvmanifolds which do not have trivial canonical bundle. Let us mention for instance the Oeljeklaus-Toma manifolds, introduced in \cite{OT}.
These complex manifolds were constructed from certain number fields, but later Kasuya showed in \cite{Ka} that they are complex solvmanifolds. There are complex solvmanifolds which admit an invariant holomorphic section of the canonical bundle, just as in the case of nilmanifolds. A classification of the Lie algebras associated to such solvmanifolds in dimension 6 is given in \cite{FOU} (although one Lie algebra is missing, see Lemma \ref{lemma: g10} below). 

	More recently, in \cite[Example 2.1]{AT} it is exhibited an example of a $4$-dimensional complex solvmanifold $(\Gamma\backslash (\R\ltimes H_3) ,J)$ with holomorphically trivial canonical bundle such that the holomorphic trivializing section is not induced by a left-invariant holomorphic $(2,0)$-form on $G:=\R\ltimes H_3$. However, this complex solvmanifold is biholomorphic to a complex nilmanifold, and this biholomorphism can be viewed as a particular case of the so-called \textit{$S$-modification} of a Lie group $G$, introduced by \cite{CFK} (see Example 2.9 of this reference). It is then natural to ask whether this is something that occurs more generally. In other words, is every complex solvmanifold with holomorphically trivial canonical bundle biholomorphic to one with an invariant trivializing holomorphic section of its canonical bundle? The main motivation for this article is to address this question. In order to achieve this, and also due to the lack of examples of this phenomenon, our goal is to construct complex solvmanifolds with holomorphically trivial canonical bundle associated to $6$-dimensional solvable Lie algebras such that the holomorphic trivializing section is not invariant.
	
	In general, we have a sufficient and necessary condition for a complex solvmanifold to admit an invariant section of their canonical bundle: according to \cite[Theorem 3.1]{AT}, a quotient of a simply connected Lie group $G$ by a uniform lattice $\Gamma$ equipped with an invariant complex structure admits an invariant holomorphic trivializing section of its canonical bundle if and only if the Koszul $1$-form $\psi$ on the Lie algebra $\g$ of $G$ vanishes identically, where $\psi$ is defined by $\psi(x)=\Tr(J \ad x)-\Tr \ad(Jx)$, for $x\in\g$.
	
	Also in terms of the Koszul $1$-form it is possible to give the following very useful algebraic obstruction (see for instance \cite[Theorem 5.3]{AT}), which also holds when $G$ is not solvable: If a complex solvmanifold has holomorphically trivial canonical bundle (or more generally, holomorphically torsion) then $\psi$ must vanish on $[\g,\g]$, or equivalently $d\psi=0$, that is, $\psi$ is closed. Notice that, according to \cite{V}, $d\psi=-2\rho$, where $\rho$ is the Chern-Ricci form associated to any left-invariant Hermitian metric $g$ on the solvable Lie group $(G,J)$. 
	 
A list of $6$-dimensional solvable Lie algebras is given in \cite{SW}. Recall that, according to \cite{Mil}, if a simply connected Lie group $G$ admits compact quotients by lattices then the Lie algebra $\g$ of $G$ must be unimodular, that is, $\Tr( \ad x)=0$ for all $x\in \g=\operatorname{Lie}(G)$. When $G$ is solvable, this condition can be strengthened. Indeed, by a result of \cite{Gar}, the Lie algebra $\g$ must be strongly unimodular, that is, $\Tr \ad x|_{\n^k / \n^{k+1}}=0$ for every $x\in \g$, $k\in \N$ where $\n$ is the nilradical of $\g$ and $\n^k$ is the $k$-th term of the descending central series of $\n$. Motivated by this we will focus on strongly unimodular Lie algebras. 
	
In Section \ref{S: strongly unimod} we first identify the six-dimensional solvable (non-nilpotent) strongly unimodular solvable Lie algebras and complete the classification of those that admit complex structures. The existence of complex structures on solvable decomposable Lie algebras of the form $\g\times \h$ with $\dim \g=\dim \h=3$ or $\dim\g=4$ and $\dim \h=2$ was established by \cite{Sroka}, while for \textit{almost nilpotent} Lie algebras (i.e. the codimension of the nilradical of the Lie algebra is 1) this was addressed on \cite{FPAA,FPAN1,FPAN2}. It follows from \cite{SW} that, in order to finish the classification, it only remains to study the case where the nilradical of the Lie algebra is $\R^4$. Since there are only a few of such Lie algebras, we can perform a case-by-case analysis, which is done in Theorem \ref{theorem: complejas}. Having this classification is also useful for the purpose of searching for Hermitian metrics on these Lie algebras, given the interest in finding Hermitian metrics that generalize the Kähler condition. 

Next we search, among the $6$-dimensional strongly unimodular Lie algebras admitting complex structures, those which admit a complex structure such that the associated Koszul $1$-form is closed. This is carried out in Theorem \ref{theorem: CRF}, where $21$ Lie algebras (in several cases, families of Lie algebras) are found, aside from those where $\psi=0$, in which cases $\psi$ is trivially closed. Furthermore, we exhibit in Lemma \ref{lemma: g10} an example of a Lie algebra with a holomorphic $(3,0)$-form which is missing from the list of \cite[Theorem 2.8]{FOU}. Interestingly enough, Theorem \ref{theorem: CRF} also furnishes a classification of six-dimensional solvable Lie groups admitting compact quotients by lattices which carry complex structures with left-invariant Chern-Ricci flat metrics.  

The goal of Section \ref{S: TCB} is to construct complex solvmanifolds with holomorphically trivial canonical bundle via non-invariant holomorphic sections. For each of the Lie algebras that admit complex structures with closed Koszul $1$-form $\psi$ we  first provide, in Theorem \ref{theorem: expression of psi}, an expression of $\psi$ in terms of the dual basis associated with the basis that defines the structure equations of the Lie algebra. This serves two purposes: first, to ensure, in the strongly unimodular case, that no other Lie algebra is missing from the list of \cite{FOU}, aside from the one found in Lemma \ref{lemma: g10}. Second, to justify working with the complex structure found in Theorem \ref{theorem: CRF}, since no new complex solvmanifolds would arise by choosing another complex structure. 

Then, we use \cite[Proposition 5.10]{AT} to construct a non-vanishing closed $(3,0)$-form $\tau$ on the associated simply connected Lie groups. In each case we are able to construct a lattice such that $\tau$ is invariant under the action of $\Gamma$, so that the canonical bundle of the corresponding complex solvmanifold is holomorphically trivial. In particular, we obtain in Theorem \ref{theorem: classification} a classification of the Lie algebras corresponding to 6-dimensional complex solvmanifolds with holomorphically trivial canonical bundle. As a consequence, we obtain a result about $6$-dimensional completely solvable complex solvmanifolds with holomorphically trivial canonical bundle (Theorem \ref{theorem: CS}). 

 In most cases we obtain examples with similar behavior to \cite[Example 2.1]{AT}, in the sense that the complex solvmanifold is biholomorphic to one with an invariant section of its canonical bundle. Nevertheless, we obtain at least one example of a complex solvmanifold that cannot be biholomorphic to a complex solvmanifold with a holomorphic invariant trivializing section (Proposition \ref{prop: S6,154}), among other interesting examples (Remark \ref{rmk: candidatos}).

The $6$-dimensional solvable (non-nilpotent) strongly unimodular Lie algebras are listed in Tables \ref{table:Decomposable 3x3}-\ref{table:nilradical codimension 2} for easier reference. The computations of Theorems \ref{theorem: complejas}, \ref{theorem: CRF}, \ref{theorem: expression of psi} and Lemma \ref{lemma: g10}, were carefully done with MAPLE 17 and can be found in full detail in \href{https://drive.google.com/drive/folders/1nggW0cmMbyYTcE53e4ifc-WtI0lx4euc?usp=sharing}{this link}.

\smallskip

\textbf{Acknowledgments.} 
The author is grateful to Adrián Andrada for his careful reading of a first version of this manuscript, to Anna Fino for useful suggestions, and to the anonymous referees for their very useful observations and suggestions that helped to improve the presentation of the paper. He would also like to thank the warm hospitality of the Dipartimento di Matematica at Università degli studi di Torino, where this work was finished. This work is supported by the PRIN
2022 project “Interactions between Geometric Structures and Function Theories” (code 2022MWPMAB).

\section{Preliminaries}\label{S: prelim}

An almost complex structure on a differentiable manifold $M$ is an automorphism
$J$ of the tangent bundle $TM$ satisfying $J^2=-\Id_{TM}$, and the existence of such a structure on $M$ forces the dimension of $M$ to be even, say $\dim_\R M=2n$. The almost complex structure $J$ is called \textit{integrable} when it satisfies the condition $N_{J}\equiv 0$, where $N_J$ is the Nijenhuis tensor given by:
\begin{equation}\label{eq:nijenhuis}
	N_{J}(X,Y) =
	[X,Y]+J([JX,Y]+[X,JY])-[JX,JY],
\end{equation}
for $X,Y$ vector fields on $M$. An integrable almost complex structure is called simply a complex structure on $M$. According to the well-known Newlander-Nirenberg theorem, a complex structure on $M$ is equivalent to the existence of a holomorphic atlas on $M$, so that $(M,J)$ can be considered as a complex manifold of complex dimension $n$. 

Given a complex manifold $(M,J)$ with $\dim_\C M=n$ its canonical bundle is defined as 
$ K_{(M,J)}=\alt^n {\mathcal T}^*_M$,  
where ${\mathcal T}^*_M$ is the holomorphic cotangent bundle of $M$. This is a holomorphic line bundle on $M$, and it is holomorphically trivial if and only if there exists a nowhere vanishing holomorphic $(n,0)$-form defined on $M$. More generally, a complex manifold $(M,J)$ is said to be \textit{holomorphically torsion} if some power $K_{(M,J)}^{\otimes k}$ is holomorphically trivial, $k\geq 1$.

Note that if $\sigma$ is a $(n,0)$-form on $M$ then $\sigma$ is holomorphic if and only if it is closed, since $d\sigma=\partial \sigma+\overline{\partial}\sigma$ and $\partial \sigma$ is a $(n+1, 0)$-form, thus $\partial \sigma=0$. 

We now move on and consider compact quotients of Lie groups by discrete subgroups. A discrete subgroup $\Gamma$ of a Lie group $G$ is called a \textit{lattice} if the quotient $\Gamma\backslash G$ has finite volume. According to \cite{Mil}, if such a lattice exists then the Lie group must be unimodular, that is, it carries a bi-invariant Haar measure. This is equivalent, when $G$ is connected, to $\Tr( \ad x)=0$ for all $x\in \g=\operatorname{Lie}(G)$ (in this case, $\g$ is called unimodular as well). When $\Gamma\backslash G$ is compact the lattice $\Gamma$ is said to be \textit{uniform}. It is well known that when $G$ is solvable then any lattice is uniform \cite[Theorem 3.1]{Rag}. 

A complex structure $J$ on a Lie group $G$ is said to be left-invariant if left translations by elements of $G$ are holomorphic maps. In this case $J$ is determined by the value at the identity of $G$. Thus, a left-invariant complex structure on $G$ amounts to a complex structure on its Lie algebra $\g$, that is, a real linear transformation $J$ of $\g$ satisfying $J^2 = -\Id$ and $N_J(x, y)=0$ for all $x, y$ in $\g$, with $N_J$ defined as in \eqref{eq:nijenhuis}. Furthermore, if $G$ is simply connected and admits a uniform lattice $\Gamma$, a left-invariant complex structure defined on $G$ induces a unique complex structure on the compact quotient $\Gamma\backslash G$ such that the projection $\pi:G\to \Gamma\backslash G$ is a local biholomorphism. Such a complex structure is called \textit{invariant}.

Let $(\Gamma\backslash G,J)$ be a $2n$-dimensional compact quotient of a simply connected Lie group $G$ by a lattice $\Gamma$, equipped with an invariant complex structure. In \cite{AT} it is characterized under which conditions the canonical bundle of $(\Gamma\backslash G,J)$ admits an \textit{invariant} trivializing holomorphic section, that is, there is a non-zero closed $(n,0)$-form defined on the Lie algebra $\g$ of $G$. According to \cite[Theorem 3.1]{AT}, this happens if and only if the Koszul $1$-form $\psi\in \g^*$ vanishes identically on $\g$, where $\psi(x)=\Tr(J\ad x)-\Tr \ad(Jx)$, for $x\in\g$. Since we will always work with unimodular Lie algebras, we may assume $\Tr \ad(Jx)=0$ for all $x\in \g$. 

While there are examples of compact complex manifolds $(\Gamma\backslash G,J)$ with holomorphically trivial canonical bundle via a holomorphic section which is not invariant, the following algebraic obstruction, given also in terms of $\psi$, turns out to be quite useful for ruling out the possibility that $(\Gamma\backslash G,J)$ has holomorphically trivial canonical bundle, via an invariant section or not:

\begin{theorem}\cite[Theorem 5.2]{AT}\label{theorem: obstruction}
If the canonical bundle of $(\Gamma\backslash G,J)$ is holomorphically trivial (or more generally, holomorphically torsion) then $\psi([\g,\g])\equiv 0$, or equivalently $d\psi=0$.
\end{theorem}

Furthermore, the canonical bundle of a compact complex quotient $\Gamma\backslash G$ cannot admit both types of trivializations, since two holomorphic trivializations must be proportional by a constant, according to \cite[Lemma 4.1]{AT}.

Assume that $G$ is simply connected and $\Gamma$ is a uniform lattice in $G$. When $G$ is solvable (resp. nilpotent) the compact quotient $\Gamma\backslash G$ is called a \textit{solvmanifold} (resp. \textit{nilmanifold}). Consequently, if $J$ is an invariant complex structure, we will call $(\Gamma\backslash G,J)$ a complex solvmanifold (resp. complex nilmanifold). These manifolds have the very nice property that $\pi_1(\Gamma\backslash G)\cong \Gamma$. 

Furthermore, the diffeomorphism class of solvmanifolds is determined by the isomorphism class of the corresponding lattices, as the following result of Mostow shows:

\begin{theorem}\cite[Theorem 3.6]{Mo}\label{solv-isom}
Let $G_1$ and $G_2$ be simply connected solvable Lie groups	with $\Gamma_i$ a lattice in $G_i$, for $i=1,2$. If $\phi:\Gamma_1\to \Gamma_2$ is an isomorphism then there exists a diffeomorphism $\tilde{\phi}: G_1\to G_2$ such that $\tilde{\phi}|_{\Gamma_1}=\phi$ and $\tilde{\phi}(\gamma g)=\phi(\gamma)\tilde{\phi}(g)$ for all $\gamma\in \Gamma_1$, $g\in G_1$. 
\end{theorem}

For a solvable Lie group the unimodularity necessary condition can be strengthened. Let $\g$ be a solvable Lie algebra and $\n$ be its nilradical. If we denote by \[ \n^0:=\n, \quad \n^1:=[\n,\n], \quad \n^{\l}:=[\n, \n^{\l-1}], \quad \l\geq 2, \] the terms of the descending central series of $\n$, the Lie algebra $\g$ is said to be \textit{strongly unimodular} if $\Tr(\ad(x)|_{\n^\l / \n^{\l+1}})=0$, for every $x\in\g$ and every $\l\in \N$. Observe that if the nilradical $\n$ of $\g$ is $r$-step nilpotent, then one has $\Tr(\ad x)=\sum_{\l=0}^{r-1} \Tr(\ad x|_{\n^\l / \n^{\l+1}}),  x\in \g.$ In particular, a strongly unimodular solvable Lie algebra is always unimodular, i.e. $\Tr \ad x=0$ for every $x\in \g$. If $\n$ is abelian, then the strongly unimodular condition is just being unimodular. By a result of \cite{Gar}, the Lie algebra of a simply connected solvable Lie group admitting compact quotients by lattices must be strongly unimodular. 


We want to construct lattices in simply connected Lie groups admitting complex structures with closed Koszul $1$-form. Since all of the Lie groups that we consider will be of the form $\R^k\ltimes_{\phi} N$, with $N$ the nilradical, we will be able to construct lattices of the form $\Gamma=\Gamma_1\ltimes_{\phi} \Gamma_2$ where $\Gamma_1\subset\R^k$ and $\Gamma_2\subset N$ are lattices of $\R^k$ and $N$ respectively. The main tool we will use is the following criterion due to \cite{Yamada}. 
\begin{theorem}\cite[Theorem 2.4]{Yamada}\label{theorem: yamada}
	Let $G=\R^k \ltimes_\phi N$ be a simply connected solvable Lie group, where $N$ is the nilradical of $G$. If there exists a rational basis $\mathcal{B}=\{X_1,\ldots,X_n\}$ of $\n$ and a basis $\{t_1,\ldots,t_k\}$ of $\R^k$ such $[d(\phi(t_j))_{1_N}]_\B$ is an integer unimodular matrix for all $1\leq j\leq k$ then $G$ has a lattice of the form $\Gamma=\text{span}_\Z\{t_1,\ldots,t_k\}\ltimes_\phi \exp^N(\text{span}_\Z \{X_1,\ldots,X_n\})$.
\end{theorem} 

Note that when the nilradical $N$ is abelian, every basis of its Lie algebra $\n$ is rational. 
When $k=1$ the Lie group $G=\R \ltimes_{\phi} N$ is called \textit{almost nilpotent}. 
Moreover, if $N$ is abelian, i.e. $N=\R^n$, then $G$ is called \textit{almost abelian}. 

In the examples we will begin with a Lie algebra $\g=\R^k\ltimes_{\varphi}\n$. To apply Theorem \ref{theorem: yamada} we need to determine the associated simply connected Lie group $G$. Let $N$ denote the simply connected nilpotent Lie group with Lie algebra $\n$. Since the exponential map $\exp:\n\to N$ is a diffeomorphism, we may assume that the underlying manifold of $N$ is $\n$ itself with the group law $x\cdot y=Z(x,y),$ where $Z(x,y)$ is the polynomial map given by the Baker-Campbell-Hausdorff formula: $\exp(x)\exp(y)=\exp(Z(x,y))$. Therefore, with this assumption, we have that $\exp:\n\to N$ is simply the identity map on $\n$ and moreover, $\Aut(\n)=\Aut(N)$.

Let $\{t_1,\ldots,t_k\}$ be a basis of $\R^k$ and denote $B_j=\varphi(t_j)\in \Der(\n)$. Then, $\exp(B_j)\in \Aut(N)$ and using \cite[Theorem 4.2]{Bock} we have that $G=\R^k\ltimes_{\phi} N$, where $\phi:\R^k\to \Aut(N)$ is the Lie group homomorphism given by \[ \phi\big(\sum_{j=1}^k  x_jt_j\big)=\exp(x_1 B_1+\cdots+x_k B_k)=\exp(x_1 B_1)\exp(x_2 B_2)\cdots \exp(x_k B_k).\] Here  $\exp$ denotes the matrix exponential after identification of $\n\cong \R^{\dim \n}$ choosing a basis of $\n$. 

Note that, in the notation of Theorem \ref{theorem: yamada}, we have that  $d(\phi(t_j))_{1_N}=\exp(B_j)=\exp(\varphi(t_j))$. Hence, in order to find lattices we need a basis $\{t_1,\ldots,t_k\}$ such that
$[\exp(\varphi(t_j))]_\B$ is an integer unimodular matrix in the rational basis $\B$ of $\n$, for all $1\leq j\leq k$.

\medskip

In order to produce examples of complex solvmanifolds with holomorphically trivial canonical bundle via a holomorphic section which is not invariant under the action of the group, we will use the following proposition, which gives a non-vanishing closed $(n,0)$-form (not necessarily left-invariant) on any simply connected solvable unimodular Lie group equipped with a left-invariant complex structure.

\begin{proposition}\cite[Proposition 5.10]{AT}\label{proposition: tau explicit}
	Let $(G,J)$ be a $2n$-dimensional simply connected solvable unimodular Lie group equipped with a left-invariant complex structure. Let $\h$ denote the kernel of $\psi:\g\to \R$ and assume that $\h\neq \g$ and $\psi([\g,\g])\equiv 0$, so that $\g=\R e_0 \ltimes \h$, and consequently $G=\R\ltimes H$, where $H$ is the unique connected normal subgroup of $G$ such that $\operatorname{Lie}(H)=\h$. Then the $(n,0)$-form \begin{equation}\label{eq: tau}
		\tau=\exp(-\tfrac{i}{2} \Tr (J \ad e_0) t) \sigma
	\end{equation} is closed, where $t$ is the coordinate of $\R$ and $\sigma$ is a left-invariant $(n,0)$-form. 
\end{proposition}

If $G$ admits a lattice $\Gamma$ such that $\tau$ is invariant under the action of $\Gamma$, we can induce a nowhere vanishing closed $(n,0)$-form on the complex manifold $(\Gamma\backslash G,J)$, and hence its canonical bundle is holomorphically trivial. 

\section{Six-dimensional solvable (non-nilpotent) strongly unimodular Lie algebras admitting a complex structure with closed Koszul form}\label{S: strongly unimod}

We begin our discussion by presenting the list of $6$-dimensional solvable (non-nilpotent) strongly unimodular Lie algebras, which is extracted from the list of all solvable Lie algebras up to dimension 6, provided in \cite{SW}, with the notation used therein. We will always identify a Lie algebra via its structure equations. For instance, the notation
$\s_{3,1}^{-1}\times \R^3=(e^{13},-e^{23},0,0,0,0)$ means that the Lie algebra $\s_{3,1}^{-1}\times \R^3$ admits a basis $\{e_1,\ldots,e_6\}$ whose dual basis $\{e^1,\ldots,e^6\}$ satisfies the following equations with respect to the Chevalley-Eilenberg differential:
$de^1=e^{13}$, $de^2 =-e^{23}$, $de^j=0\; (3\leq j\leq 6).$ 

We first list the unimodular Lie algebras and from them, we keep the strongly unimodular ones. The analysis can be carried out through a straightforward verification of the definitions of unimodularity and strongly unimodularity for each Lie algebra.

\subsection{Six-dimensional solvable (non-nilpotent) strongly unimodular Lie algebras}
 
Six-dimensional solvable unimodular Lie algebras can be divided into two groups, according to whether they are decomposable or indecomposable. 
	
	\subsection*{Decomposable Lie algebras}
	In the decomposable case, $\g\times \h$ is solvable, unimodular, or strongly unimodular if and only if both $\g$ and $\h$ are, respectively, solvable, unimodular, or strongly unimodular. Moreover, for $\g\times \h$ to be non-nilpotent, at least one of the two factors must not be nilpotent. In dimension $1$, the only Lie algebra is the abelian algebra $\R$ and in dimension $2$, the only unimodular Lie algebra is $\R^2$. Thus, we can further subdivide the decomposable case into three subcases:
	\begin{itemize}
		\item $\g\times \h$ where $\g$ and $\h$ are $3$-dimensional strongly unimodular solvable Lie algebras (with at least one of them being non-nilpotent), or
		\item $\g\times \R^2$ where $\g$ is a $4$-dimensional indecomposable strongly unimodular solvable non-nilpotent Lie algebra, or
		\item $\g\times \R$ where $\g$ is a $5$-dimensional indecomposable strongly unimodular solvable non-nilpotent Lie algebra.
	\end{itemize}
	
	In dimension $3$, the solvable unimodular Lie algebras are $\R^3$, $\n_{3,1}$, $\s_{3,1}^{-1}$ and $\s_{3,3}^0$. Since all of these have abelian nilradical and are unimodular, they are strongly unimodular. The Lie algebra $\n_{3,1}$ is the $3$-dimensional Heisenberg Lie algebra spanned by $\{e_1,e_2,e_3\}$ with Lie bracket defined by $[e_2, e_3]=e_1$. From now on we will denote $\h_3:=\n_{3,1}$.
	
	In dimension $4$, the solvable non-nilpotent unimodular Lie algebras are: $\s_{4,3}^{a,-(1+a)}$  {\tiny$(a\in (-1,\tfrac12])$}, \; $\s_{4,4}^{-2}, \; \s_{4,5}^{\alpha, -\frac{\alpha}{2}} $ {\tiny($\alpha>0$)}$, \;\s_{4,6}, \; \s_{4,7},$ and they are all strongly unimodular. 
	
	In dimension $5$, the solvable non-nilpotent unimodular Lie algebras can be divided according to their nilradical:
	\begin{itemize} 
		\item $\nil(\g)=\R^4$: $\s_{5,3}^{-1},\; \s_{5,4}^0,\; \s_{5,6}^{-1},\; \s_{5,7}^{-3},\; \s_{5,8}^0,\; \s_{5,9}^{a,b,-(a+b+1)}$  {\tiny($0<|a+b+1|\leq |b|\leq |a|\leq 1)$},\, $\s_{5,10}^{a,-(a+2)}$ {\tiny($a\leq -1, a\neq -2$)},$ \;
		\s_{5,11}^{\alpha, \beta, -\frac{\alpha+\beta}{2}}$ {\tiny($\alpha>0,\, \beta\neq 0,\, |\beta|\leq \alpha$)}$, \; \s_{5,12}^{-1,\beta}$ {\tiny($\beta>0$)}$, \; \s_{5,13}^{\alpha,-\alpha,\gamma} $ {\tiny($0<\gamma\leq 1,\, 0\leq \alpha$)}. Since all of these have abelian nilradical and are unimodular, they are automatically strongly unimodular.
		
		\item $\nil(\g)=\h_3\times \R$: $\s_{5,15}, \; \s_{5,16}, \; \s_{5,22}^{a,-(a+2)}$ {\tiny ($a\in (-1,1]\setminus\{0\}$)}$, \; \s_{5,23}^{-\frac32}, \; \s_{5,24}^{-4}, \; \s_{5,25}^{\alpha,-4\alpha}$ {\tiny ($\alpha<0$)}$, \; \s_{5,30}^{-2}$. Among all of these, only $\s_{5,15}, \s_{5,16}$ are strongly unimodular.
		
		\item $\nil(\g)=\n_{4,1}$: $\s_{5,35}^{-\frac43}$, which is not strongly unimodular.
		
		\item $\nil(\g)=\R^3$: $\s_{5,41}^{-1,-1}, \s_{5,43}^{-2,0}$, which are strongly unimodular as they are already unimodular and have an abelian nilradical.
	\end{itemize} 
	
\subsection*{Indecomposable Lie algebras} 
The unimodular $6$-dimensional solvable (non-nilpotent) Lie algebras can be divided into several groups based on their nilradical.

\begin{itemize}
	\item $\mathfrak{nil}(\g)=\R^5$: $\s_{6,4}^{-1},\;\s_{6,5}^0, \; \s_{6,7}^{-\frac12}, \; \s_{6,8}^{a,-(a+1)}$ {\tiny($-1<a\leq -\tfrac12$)}$, \; \s_{6,9}^{\alpha,-\frac{\alpha}{2}} ${\tiny($\alpha>0$)}$, \; \s_{6,11}^{-\frac32}, \; \s_{6,12}^{-\frac14},\\
	\s_{6,13}^{a,-(3a+1)}$ {\tiny($a\in [-\tfrac23, 0)\setminus\{-\tfrac13\}$)},$\;\s_{6,14}^{a,-(a+\frac12)}$ {\tiny($a\leq -\tfrac14, \, a\neq -\tfrac12$)}$,\; \s_{6,15}^{\alpha, -\frac{3\alpha}{2}}${\tiny($\alpha>0$)}$,\; \s_{6,16}^{\alpha, -4\alpha}$ {\tiny($\alpha<0$)}$,\\ \s_{6,17}^{a,b,c,-(a+b+c+1)}${\tiny ($0<|a+b+c+1|\leq |c|\leq |b|\leq |a|\leq 1$)}$,\;  \s_{6,18}^{a,b,-(2a+b+1)} ${\tiny($0<|2a+b+1|\leq |b|\leq 1,\, a\neq 0$)}$,\\
	\s_{6,19}^{\alpha,\beta,\gamma,-\frac{\alpha+\beta+\gamma}{2}}${\tiny($0<|\gamma|\leq|\beta|\leq \alpha$)}$,\; \s_{6,20}^{\alpha,\beta,-\frac{2\alpha+\beta}{2}} ${\tiny($0<\alpha,\, \beta \neq 0$)}$, \\ \s_{6,21}^{\alpha,\beta,\gamma,-2(\alpha+\beta)} ${\tiny($\alpha<-\beta, \, \gamma\in (0,1]. \text{ If } \gamma=1 \text{ then } \alpha\leq \beta$)}. Since all of these have an abelian nilradical and are already unimodular, they are automatically strongly unimodular. 
	
	\item $\mathfrak{nil}(\g)=\h_3\times \R^2$:   
	$\s_{6,24}, \; \s_{6,25},\; \s_{6,30},\;\s_{6,31},\; \s_{6,32}^{-1},\; \s_{6,34}^0, \; \s_{6,43}, \; \s_{6,44}, \;
	\s_{6,45}^{a,-1}$ {\tiny  ($a\neq 0$)}$,\\
	\s_{6,46}^{\alpha,-\alpha}$ {\tiny ($\alpha>0$)}$, \, \s_{6,47}^{-1}, \; \s_{6,51}^{\alpha,0}$ {\tiny ($\alpha>0$)}$, \; \s_{6,52}^{0,\beta}$ {\tiny($\beta>0$)}$, \; \s_{6,54}^{-\frac12}, \; \s_{6,57}^{-\frac15},  \\ \s_{6,66}^{a,b,-(2a+2b+1)}${\tiny ($0<|b|\leq |a|,\, 0<|2a+2b+1|\leq 1. \text{ If } |a|=|b| \text{ then } a\neq -b$)}$, \,
	\s_{6,67}^{a,-(2a+3)}$ {\tiny($a\neq -1,0,-\tfrac32$)}$, \\  \s_{6,68}^{a,-(4a+1)} ${\tiny ($a\in [-\tfrac12, 0)\setminus\{-\tfrac14\}$)}$, \;
	\s_{6,69}^{\alpha,\beta,-(4\alpha+\beta)} $ {\tiny($\alpha\neq 0,\, 0<|4\alpha+\beta|\leq |\beta|. \text{ If } |4\alpha+\beta|=|\beta| \text{ then } 2\alpha\geq -\beta$)}$, \;
	\s_{6,70}^{-2}, \\
	\s_{6,71}^{a,-(a+\frac13)}$ {\tiny ($a\leq -\tfrac16,\, a\neq -\tfrac13$)}$,  \; \s_{6,72}^{-\frac43},  \;  \s_{6,73}^{-\frac16},  \;   \s_{6,74}^{a,-(a+1)} $ {\tiny($0<|a|\leq 1,\, a\neq -1$)}$,  \; \s_{6,75}^{-2},  \;  \s_{6,76}^{\alpha,-6\alpha} ${\tiny($\alpha>0$)}$, \\ \s_{6,77}^{\alpha,-2\alpha} $ {\tiny($\alpha>0$)}$,\; \s_{6,78}^{\alpha,\beta,-(\alpha+\beta)} $ {\tiny($0<|\beta|\leq |\alpha|, \, \beta\neq -\alpha, \, \alpha+\beta<0$)}$, \quad\s_{6,79}^{\alpha,-2\alpha}$ {\tiny ($\alpha<0$)}$,\; \s_{6,80}^{\alpha,-2\alpha,\gamma}$ {\tiny ($\alpha,\gamma>0$)}$,  \; \s_{6,81}^{-3},  \\ \s_{6,86}^{a,-(2a+1)}$ {\tiny($a\in [-1,0)\setminus\{-\tfrac12\}]$)}$,  \; \s_{6,87}^{-\frac13},\;
	\s_{6,88}^{-3}, \; \s_{6,89}^{-1}, \; \s_{6,90}^{\alpha, -\alpha}$ {\tiny ($\alpha<0$)}, 
	of which only the Lie algebras $\s_{6,24},\; \s_{6,25},\; \s_{6,30}, \;\s_{6,31}, \;\s_{6,32}^{-1},\;
	\s_{6,34}^0, \;\s_{6,43}, \s_{6,44},\;\s_{6,45}^{a,-1}, \; \s_{6,46}^{\alpha,-\alpha}$ {\tiny ($\alpha>0$)}$ , \; \s_{6,47}^{-1}, \; \s_{6,51}^{\alpha,0}$ {\tiny($\alpha>0$)}$, \\ \s_{6,52}^{0,\beta}$ {\tiny($\beta>0$)}, are strongly unimodular. 
	
	\item $\nil(\g)=\n_{4,1}\times \R$: $
	\s_{6,96}^{2}, \; \s_{6,102}^{-7}, \; \s_{6,105}^{a,-(3a+4)} ${\tiny($a\neq 0, -2, -\tfrac43$)}$, \;\s_{6,106}^{-\frac{5}{3}}, \, \s_{6,107}^{-1}, \;
	\s_{6,108}^{-\frac32}, \;\s_{6,110}^{-4}, \;\s_{6,112}^{-3}, \\ \s_{6,116}^{-3},$ 
	of which none is strongly unimodular. 
	\item $\nil(\g)=\n_{5,1}$:  $\s_{6,124}^2, \; \s_{6,128}^{-\frac32}, \; \s_{6,131}^{a,2+2a}$ {\tiny($0<|a|\leq1,\, a\neq -1,-\tfrac12$)}$, \; \s_{6,132}^{-3}, \; \s_{6,133}^4, \; \s_{6,134}^{\alpha,4\alpha}$ {\tiny($\alpha>0$)}, \; $\s_{6,139}^{-\frac12}, \; \s_{6,140}^{-1}, \; \s_{6,145}^{0}, \; \s_{6,146}^{-1}, \; \s_{6,147}^{0}$,
	of which just $\s_{6,140}^{-1}, \s_{6,145}^{0}, \s_{6,146}^{-1}, \s_{6,147}^{0}$ are strongly unimodular.
	
	\item $\nil(\g)=\n_{5,2}$: $\s_{6,151}$, \; $\s_{6,152}^{\pm 1}, \; \s_{6,154}^0, \; \s_{6,155}^{-1}$, which are all strongly unimodular. Note that
	$\s_{6,152}^{-1} \simeq \s_{6,152}^1$ (\cite[Remark 2.1]{FPAN2}).
	
	\item$\nil(\g)=\n_{5,3}$:\quad $
	\s_{6,158}, \; \s_{6,159}, \; \s_{6,160}, \; \s_{6,161}^{\pm 1}, \; \s_{6,162}^a $ {\tiny($0<|a|\leq 1$)}$, \; \s_{6,163},\; \s_{6,164}^{\alpha}  $ {\tiny($\alpha>0$)}$, \; \s_{6,165}^{\alpha}                  $ {\tiny($\alpha>0$)}$, \;  \s_{6,166}^{\alpha}  $ {\tiny($0<|\alpha|\leq 1$)}$, \;\s_{6,167}$,
	which are all strongly unimodular.
	
	\item $\nil(\g)=\n_{5,4}$: $\s_{6,183}$, which is not strongly unimodular.
	
	\item $\nil(\g)=\n_{5,5}$: $\s_{6,192}^{-\frac74}$, which is not strongly unimodular.
	
	\item $\nil(\g)=\R^4$: $\s_{6,204}^{-1,-1},  \;  \s_{6,208}^{0,-2},  \;  \s_{6,213}^{a,-(a+1),c,-(c+1)} $ {\tiny($a\leq-\tfrac12$)}$, \; 
	\s_{6,214}^{a,-(a+2),-1} $ {\tiny($a\in \R$)}$, \;\s_{6,215}^{-2,-1}, \; \s_{6,216}^{\alpha, -2\alpha, -1} $ {\tiny($\alpha<0$)}$, \;
	\s_{6,217}^{\alpha,-\alpha,\gamma,-(\gamma+2)} $ {\tiny($\alpha>0 \text{ or } (\alpha=0 \text{ and }\gamma \geq -1,\, \gamma \neq 0)$)}$  , \;\s_{6,224}^{\alpha,-\alpha,-1} $ {\tiny($\alpha\geq 0$)}$, \;
	\s_{6,226}^{0,\beta, -1} $ {\tiny($0<\beta\leq 1$)}$, \; \s_{6,227}^{0,-1}, \\ \s_{6,228}^{-\beta,\beta,\gamma,-\gamma} $ {\tiny($0\leq \beta\leq \gamma,\, 0<\gamma$)}, which, being unimodular and having an abelian nilradical, are automatically strongly unimodular.
	\item $\nil(\g)=\h_3 \times \R$: $\s_{6,234}^{-2,-2}, \; \s_{6,239}^{-4,0}$, which are not strongly unimodular.
\end{itemize}

As a consequence of this analysis, we have identified the $6$-dimensional solvable non-nilpotent strongly unimodular Lie algebras, which we present along with their structure equations in \S\ref{S: algebras}. 

The decomposable case is detailed in Tables \ref{table:Decomposable 3x3}, \ref{table:Decomposable 4+R2}, \ref{table:Almost abelian 5+R} and \ref{table:Decomposable 5+R}, while the indecomposable case is found in Tables \ref{table:Almost abelian indecomposable}, \ref{table:almost nilpotent} and \ref{table:nilradical codimension 2}. The nilradicals appearing in Table \ref{table:almost nilpotent} labeled as $\n_{5,1}$, $\n_{5,2}$ and $\n_{5,3}$, have the following structure equations:
\begin{align*} 
	\n_{5,1}&=(-e^{35},-e^{45},0,0,0), \\
	\n_{5,2}&=(-e^{35},-e^{34},-e^{45},0,0), \\
	\n_{5,3}&=(-e^{24}-e^{35},0,0,0,0).
	\end{align*}
 Notice that the only nilradical of codimension 2 that occurs is $\R^4$, and that there are no strongly unimodular Lie algebras whose nilradical has codimension $>2$.

\subsection{6-dimensional strongly unimodular Lie algebras admitting a complex structure}\label{S: complex}

We complete next the classification of $6$-dimensional solvable (non-nilpotent) strongly unimodular Lie algebras that admit a complex structure. As we mentioned in the introduction, the decomposable $3\times 3$ and $4\times 2$ cases and the almost nilpotent case are already done. Taking into account our analysis of the strongly unimodular Lie algebras, all the remaining Lie algebras we have to examine have nilradical $\R^4$. We obtain the following theorem.

\begin{theorem}\label{theorem: complejas}
	A six-dimensional solvable (non-nilpotent) strongly unimodular Lie algebra $\g$ with nilradical $\R^4$ admits a complex structure if and only if $\g$ is isomorphic to one among:
	\begin{align*} 
		&\s_{3,3}^0 \times \s_{3,3}^0, \quad \s_{3,1}^{-1}\times \s_{3,3}^0, \quad \s_{5,43}^{-2,0}\times \R, \quad  \s_{6,213}^{a,-(a+1),c,-(c+1)} \,(\{a,c\}\in \{ \{-\tfrac12,-\tfrac12\}\cup \{-1,-2\}\}), \\
		&\s_{6,216}^{\alpha,-2\alpha,-1} (\alpha<0), \quad  \s_{6,217}^{\alpha,-\alpha,\gamma,-(\gamma+2)} (\alpha>0 \text{ or } (\alpha,\gamma)=(0,-1)),\quad
		\s_{6,224}^{\alpha,-\alpha,-1} (\alpha>0), \\
		&\s_{6,226}^{0,\beta,-1} (0<\beta\leq 1), \quad \s_{6,227}^{0,-1}, \quad \s_{6,228}^{-\beta,\beta,\gamma,-\gamma} (0\leq \beta\leq \gamma, 0<\gamma).  \end{align*}
	\end{theorem}
\begin{proof} 
	First we provide in the following table an explicit example of a complex structure for each Lie algebra listed in the statement. Complex structures for $\s_{3,3}^0\times \s_{3,3}^0$ and $\s_{3,1}^{-1}\times \s_{3,3}^0$ already appeared in \cite{Sroka} and a complex structure for $\s_{5,43}^{-2,0}\times \R$ was recently given in \cite{BF}.
	
	\begin{table}[H]
			\renewcommand{\arraystretch}{1.3}
\small{	\begin{tabular}{|l|l|} \hline 
			$\g$  & Complex structure \\ \hline 
			$\s_{3,3}^0 \times \s_{3,3}^0$ & $Je_1=e_2, \; Je_3=e_6, \; Je_4=e_5$ \\ \hline 
			$\s_{3,1}^{-1}\times \s_{3,3}^0$ & $Je_1=e_2+e_4, \; Je_2=e_1+e_5, \; Je_3=e_6$ \\ \hline
			$\s_{5,43}^{-2,0}\times \R$ & $Je_1=e_4, \; Je_2=e_3, \; Je_5=e_6$ \\ \hline & \\ [-1.3em] 	$\s_{6,213}^{-\frac12,-\frac12,-\frac12,-\frac12}$ & $Je_1=e_2, \; Je_3=e_6, \; Je_4=e_5$ \\ 
			$\s_{6,213}^{-1,0,-2,1}$ & $Je_1=e_6, \; Je_2=e_3, \; Je_4=-2e_5+e_6$ \\
			$\s_{6,213}^{-2,1,-1,0}$ & $Je_1=e_5, \; Je_2=-e_4, \; Je_3=e_5-2e_6$\\ \hline 
			$\s_{6,216}^{\alpha,-2\alpha,-1}$ $(\alpha<0)$ & $Je_1=e_2, \; Je_3=e_5-2\alpha e_6, \; Je_4=e_5 $ \\ \hline 
			$\s_{6,217}^{\alpha,-\alpha,\gamma,-(\gamma+2)}$ $(\alpha>0)$ & $Je_1=e_2, \; Je_3=-(\gamma+2)e_5+\alpha e_6, \; Je_4=-\gamma e_5+\alpha e_6$ \\
			$\s_{6,217}^{0,0,-1,-1}$ & $Je_1=e_2, \; Je_3=e_4, \; Je_5=e_6$ \\\hline
			$\s_{6,224}^{\alpha,-\alpha,-1}$ ($\alpha>0$) & $Je_1=e_2, \; Je_3=\frac{1}{\alpha}e_5+e_6, \; Je_4=e_6$ \\ \hline 
			$\s_{6,226}^{0,\beta,-1}$ ($0<\beta\leq1$) & $Je_1=e_2, \; Je_3=e_4, \; Je_5=e_6$ \\ \hline 
			$\s_{6,227}^{0,-1}$ & $Je_1=e_2, \; Je_3=e_5, \; Je_4=-e_6$ \\  \hline
			 $\s_{6,228}^{-\beta,\beta,\gamma,-\gamma}$ ($0\leq\beta\leq\gamma, 0<\gamma$) & $Je_1=e_2, \; Je_3=e_4, \; Je_5=e_6$ \\ \hline
		\end{tabular}}
	\end{table}
	
	For the remaining Lie algebras with nilradical $\R^4$, we prove that none of them admit complex structures by performing explicit computations. We consider a generic endomorphism 
	
	{\tiny\begin{equation}\label{eq: J}  J=\matriz{a_1&a_7&a_{13}&a_{19}&a_{25}&a_{31}\\ a_2 &a_8 & a_{14} &a_{20} &a_{26} &a_{32}\\ a_3 &a_9 & a_{15} & a_{21} & a_{27} & a_{33} \\ a_4 & a_{10} & a_{16} & a_{22} & a_{28} & a_{34} \\ a_5 & a_{11} & a_{17} & a_{23} & a_{29} & a_{35} \\ a_6 & a_{12} & a_{18} & a_{24} & a_{30} & a_{36}},
		\end{equation}} written as a matrix with respect to the basis $\{e_1,\ldots,e_6\}$ which defines the structure equations. We then impose the conditions $J^2=-\Id$ and $N_J\equiv 0$ and show that they lead to a contradiction. Specifically, we set $N_{ijk}:=e^k(N_J(e_i,e_j))$ and $B_{ij}:=(J^2+\Id)_{ij}$ and show that there is some expression in terms of $N_{ijk}$ and $B_{ij}$, for some $i,j,k$, that cannot vanish. In the decomposable case, in \cite{Sroka} it was proved that $\s_{3,1}^{-1}\times \s_{3,1}^{-1}$ does not admit complex structures. The only remaining decomposable algebra to be analyzed is $\s_{5,41}^{-1,-1}\times \R=(e^{14},e^{25},-e^{34}-e^{35},0,0,0)$, which we examine next. We start by noting that $a_{34}=0$. This follows from the vanishing of
	\begin{align*} 
		N_{164}&=a_4 a_{34},& \quad N_{165}&=a_5 a_{34},\\
		N_{144}&=(a_1+a_{22})a_4-a_3a_{16},& \quad N_{145}&=(a_1+a_{22})a_5-a_3a_{17},\\ 
		N_{141}&=a_1^2-a_3 a_{13}+a_4a_{19}+1,& \quad N_{353}&=-(a_{16}+a_{17})a_{27}+a_9a_{14}-a_{15}^2-1,\\
		N_{361}&=(-2a_{34}-a_{35})a_{13}+a_{16}a_{31},& \quad N_{362}&=(-a_{34}-2a_{35})a_{14}+a_{17}a_{32}.
	\end{align*}
	Substituting $a_{34}=0$ in all $N_{ijk}$, we then deduce that $a_{35}=0$, by looking at:
	{\small\begin{align*}
		N_{265}&=a_{11}a_{35},\quad 
		N_{365}=-a_{17}a_{35}, \quad
		N_{362}=-2a_{14}a_{35}+a_{17}a_{32},\quad 
		N_{252}=a_{11}a_{26}-a_9a_{14}+a_8^2+1.
	\end{align*}}
	Next, we compute \begin{align*} 
		N_{466}+N_{566}+B_{66}&=[a_{18}a_{33}-a_{31}a_6]+[-a_{12}a_{32}+a_{18}a_{33}]+[a_{36}^2+a_{31}a_6+a_{12}a_{32}+a_{18}a_{33}]\\
		&=3a_{18}a_{33}+a_{36}^2+1.
	\end{align*} 
	The vanishing of this expression implies $a_{18}\neq 0$ and $a_{31}a_{32}\neq 0$ (the latter follows from $0=N_{466}=N_{566}$). At this point, we obtain $a_4=a_5=a_{16}=a_{17}=0$ due to the vanishing of 
	\[ N_{161}=a_4 a_{31}, \quad N_{162}=a_5 a_{32}, \quad N_{361}=a_{16}a_{31}, \quad N_{362}=a_{17}a_{32},\] 
	From the vanishing of 
	$N_{141}=a_1^2-a_3a_{13}+1$ and  $N_{151}=-a_3a_{13}+a_2a_7,$ we deduce that
	$a_2a_3 a_{13}\neq 0$. Then, by adding and substracting  \[\tfrac{N_{143}}{a_3}=a_1-a_{15}+2a_{22}+a_{23}\quad \text{and}\quad \tfrac{N_{341}}{a_{13}}=a_1-a_{15}-2a_{22}-a_{23},\] we find that $a_{15}=a_1$ and $a_{23}=-2a_{22}$. The contradiction arises since \begin{align*} 
			-(a_1-3a_{22})N_{142}+a_3 N_{342}+a_2 N_{141}&=-(a_1-3a_{22})[(a_1+3a_{22})a_2-a_3a_{14}]\\
			 &\quad +a_3[(-a_1+3a_{22})a_{14}+a_2a_{13}]+a_2 N_{141}\\
			&=9a_2(a_{22}^2+1) \neq 0.
			\end{align*}
	
	In the indecomposable case one has to perform the same analysis, for the remaining Lie algebras in Table \ref{table:nilradical codimension 2}. The computations justifying that the Lie algebras $\s_{6,204}^{-1,-1}$, $\s_{6,208}^{0,-2}$, $\s_{6,213}^{a,-(a+1),c,-(c+1)}$ ($(a,c)\neq(-\tfrac12,-\tfrac12),(-1,-2),(-2,-1))$,  $\s_{6,214}^{a,-(a+2),-1}$, $\s_{6,215}^{-2,-1}$, $\s_{6,217}^{0,0,\gamma,-(\gamma+2)}$ ($\gamma>-1$), $\s_{6,224}^{0,0,-1}$ do not admit complex structures can be found in the MAPLE files.
\end{proof} 

\subsection{Complex structures with closed Koszul form}\label{S: CRF}
We recall that, due to Theorem \ref{theorem: obstruction}, a necessary condition for a complex solvmanifold to have holomorphically trivial canonical bundle is that the associated Koszul $1$-form defined on the Lie algebra is closed. In this subsection, we classify all six-dimensional solvable (non-nilpotent) strongly unimodular Lie algebras admitting complex structures satisfying this condition. First, we will separate several Lie algebras which are known to have complex structures with $\psi \equiv 0$, or equivalently, a non-zero closed $(3,0)$-form, as in this case $\psi$ is trivially closed. In \cite[Theorem 2.8]{FOU}, the authors provide a list of these Lie algebras. They obtain nine Lie algebras (in fact one of them is a monoparametric family): $\g_1, \g_2\, (\alpha\in \R_{\geq 0}), \ldots, \g_9$. However, as we show next, there is at least one Lie algebra that is missing from this list. Our subsequent calculations will confirm that this is the only missing Lie algebra, in the strongly unimodular case.

\begin{lemma}\label{lemma: g10}
	The Lie algebra $\g_{10}:=\s_{6,147}^0$ admits a complex structure $J$ such that its associated Koszul $1$-form vanishes identically on $\g_{10}$.
\end{lemma}
\begin{proof}
	It is straightforward to verify that the almost complex structure $J$ defined on $\g_{10}$ by $Je_1= e_2$, $Je_3=-e_4$, and $Je_5=-2e_6$ is integrable and satisfies $\psi=0$.
\end{proof}

\begin{remark}
	It was observed in \cite[Remark 2.3]{FPAN2} that the simply connected Lie group associated to $\g_{10}$ admits lattices, hence it gives rise to complex solvmanifolds with holomorphically trivial canonical bundle. Moreover, $\g_{10}$ is not isomorphic to any of the Lie algebras $\g_1,\ldots,\g_9$ since $b_2(\g_{10})=2$, while none of the other algebras exhibit this property. 
\end{remark}

Next, we analyze the remaining Lie algebras that admit complex structures.

\begin{theorem}\label{theorem: CRF}
	Let $\g$ be a $6$-dimensional solvable (non-nilpotent) strongly unimodular Lie algebra different from $\g_i, 1\leq i\leq 10$. Then $\g$ admits a complex structure $J$ such that the associated Koszul $1$-form $\psi$ vanishes on $[\g,\g]$ if and only if $\g$ is isomorphic to one of the following:
\begin{align*}
	&\s_{3,3}^{0}\times \R^3, \quad \s_{3,3}^{0} \times \h_3, \quad \s_{3,3}^{0} \times \s_{3,3}^{0}, \quad 	\s_{4,7}\times \R^2, \quad \s_{5,4}^{0} \times \R, \quad \s_{5,8}^{0} \times \R, \quad \s_{5,11}^{\alpha,\alpha,-\alpha} \times \R \; (\alpha>0),\\ 
	&\s_{5,13}^{\alpha,-\alpha,\gamma}\times \R \; (0<\gamma<1, 0\leq \alpha), \quad \s_{5,16}\times \R, \quad \s_{6,25}, \quad \s_{6,44}, \quad \s_{6,52}^{0,\beta} \; (\beta>0, \beta \neq 1),\\ 
	&\s_{6,145}^{0}, \quad \s_{6,154}^{0}, \quad \s_{6,159}, \quad \s_{6,165}^{\alpha} \;(\alpha>0), \quad \s_{6,166}^{\alpha} \; (0<|\alpha|<1), \quad \s_{6,167}, \quad \s_{6,217}^{0,0,-1,1},\\ &\s_{6,226}^{0,\beta,-1} \; (0<\beta<1), \quad \s_{6,228}^{-\beta,\beta,\gamma,-\gamma} \; (0\leq \beta\leq \gamma, 0<\gamma). 
	\end{align*}
\end{theorem}

\begin{proof} 

We begin by exhibiting in the following table a complex structure such that $\psi([\g,\g])\equiv 0$ for each Lie algebra $\g$ of the list in the statement.

	\begin{table}[H]
	\renewcommand{\arraystretch}{1.1}
\small{	\begin{tabular}{|l|l|} \hline 
		$\g$  & Complex structure with $d\psi=0$\\ \hline 
		$\s_{3,3}^{0} \times \R^3$ & $Je_1=e_2, \; Je_3=e_4, \; Je_5=e_6$ \\ \hline 
		$\s_{3,3}^0 \times \h_3$ & $Je_1=e_2, \;  Je_3=e_4, \;  Je_5=e_6$ \\ \hline 
		$\s_{3,3}^{0} \times \s_{3,3}^0$ & $Je_1=e_2, \;  Je_3=e_6, \;  Je_4=e_5$ \\ \hline 
		$\s_{4,7}\times \R^2$ & $Je_1=e_4, \;  Je_2=e_3, \;  Je_5=e_6$ \\ \hline 
		$\s_{5,4}^0 \times \R$ & $Je_1=e_6,  \; Je_2=e_5, \;  Je_3=e_4$ \\ \hline 
		$\s_{5,8}^0 \times \R$ & $Je_1=e_2, \;  Je_3=e_4, \;  Je_5=e_6$ \\ \hline 
		$\s_{5,11}^{\alpha,\alpha,-\alpha} \times \R$ ($\alpha>0$) & $Je_1=e_2, \;  Je_3=e_4, \;  Je_5=e_6$ \\ \hline 
		$\s_{5,13}^{\alpha,-\alpha,\gamma}\times \R$ ($0<\gamma<1, 0\leq \alpha$) & $Je_1=-e_2, \;  Je_3=e_4,  \; Je_5=e_6$ \\ \hline 
		$\s_{5,16}\times \R$ & $Je_1=e_6, \;  Je_2=e_3, \;  Je_4=e_5$ \\ \hline 
		$\s_{6,25}$ & $Je_1=e_5, \;  Je_2=e_3, \;  Je_4=e_6$ \\ \hline 
		$\s_{6,44}$ & $Je_1=e_6, \;  Je_2=e_3, \;  Je_4=e_5$ \\ \hline 
		$\s_{6,52}^{0,\beta}$ ($\beta>0, \beta\neq 1$) & $Je_1=e_6, \; Je_2=e_3, \; Je_4=e_5$ \\ \hline 
		$\s_{6,145}^0$ & $Je_1=e_2, \;  Je_3=e_4,  \; Je_5=e_6$ \\ \hline 
		$\s_{6,154}^0$ & $Je_1=-e_2, \;  Je_3=e_6,  \; Je_4=e_5$ \\ \hline 
		$\s_{6,159}$ & $Je_1=e_6, \;  Je_2=e_4,  \; Je_3=e_5$ \\ \hline 
		$\s_{6,165}^{\alpha}$ ($\alpha>0$) & $Je_1=e_6, \;  Je_2=e_3, \;  Je_4=e_5$ \\ \hline 
		$\s_{6,166}^{\alpha}$ ($0<|\alpha|<1$) & $Je_1=e_6,  \; Je_2=e_4, \;  Je_3=e_5$ \\ \hline 
		$\s_{6,167}$ & $Je_1=e_6, \;  Je_2=e_3,  \; Je_4=e_5$ \\ \hline 
		$\s_{6,217}^{0,0,-1,-1}$ & $Je_1=e_2,  \; Je_3=e_4,  \; Je_5=e_6$ \\ \hline 
		$\s_{6,226}^{0,\beta,-1}$ ($0<\beta<1$) & $Je_1=-e_2,  \; Je_3=e_4, \;  Je_5=e_6$ \\ \hline 
		$\s_{6,228}^{-\beta,\beta,\gamma,-\gamma}$ ($0\leq\beta\leq\gamma$, $0<\gamma$) & $Je_1=e_2, \;  Je_3=e_4, \;  Je_5=e_6$ \\ \hline
	\end{tabular}}
\end{table}

Next, we show for the remaining Lie algebras that there is no complex structure with closed Koszul form. We achieve this by considering a generic endomorphism $J$ as in \eqref{eq: J} written as a matrix with respect to the basis $\{e_1,\ldots,e_6\}$ which defines the structure equations. Then we impose the conditions $J^2=-\Id$, $N_J\equiv 0$, and $\psi([\g,\g])\equiv 0$, and show that they lead to some contradiction. We define $N_{ijk}$ and $B_{ij}$ as we did in the proof of Theorem \ref{theorem: complejas}. 
For instance, for the Lie algebra
$\s_{4,3}^{-\frac12,-\frac12}\times \R^2=(e^{14},-\tfrac12 e^{24}, -\tfrac12 e^{34},0,0,0),$
	we obtain $\psi=-a_4 e^1+\tfrac{a_{10}}{2}e^2+\tfrac{a_{16}}{2}e^3+(a_1-\tfrac{a_8+a_{15}}{2})e^4.$ Thus, $\psi([\g,\g])\equiv 0$ if and only if $a_4=a_{10}=a_{16}=0$. After substitution, the following expressions do not vanish simultaneously, leading to a contradiction:
	\begin{align*}
		N_{152}&=\tfrac32 a_2 a_{28},\quad N_{153}=\tfrac32 a_3 a_{28},\quad&
		N_{162}&=\tfrac32 a_2 a_{34}, \quad N_{163}=\tfrac32 a_3 a_{34},\\
		N_{141}&= a_1^2-\tfrac12 a_2 a_7-\tfrac12 e_3 a_{13}+1, \quad &
		B_{44}&=a_{22}^2+a_{23}a_{28}+a_{24}a_{34}+1.
	\end{align*} 
	The analysis to show that the Lie algebras $\s_{4,5}^{\alpha,-\frac12 \alpha} \times \R^2$, $\s_{4,6}\times \R^2$, $\s_{5,43}^{-2,0} \times \R$, $\s_{6,14}^{-\frac14, -\frac14}$, $\s_{6,16}^{\alpha,-4\alpha}$, $\s_{6,17}^{1,b,b,-2(b+1)}$, $\s_{6,18}^{1,-\frac32,-\frac32}$, $\s_{6,19}^{\alpha,\alpha,\gamma,-\alpha-\frac{\gamma}{2}}$, $\s_{6,20}^{\alpha,\alpha,-\frac{3}{2}\alpha}$,  $\s_{6,21}^{\alpha,\beta,\gamma,-2(\alpha+\beta)}$, $\s_{6,51}^{0,\alpha}$, $\s_{6,158}$, $\s_{6,164}^{\alpha}$,  $\s_{6,213}^{a,-a-1,c,-c-1}$ (where {\tiny$(a,c)=(-\frac12,-\frac12), (-1,-2)$, or $(-2,-1)$}), $\s_{6,216}^{\alpha,-2\alpha,-1}$, $\s_{6,217}^{\alpha,-\alpha,\gamma,-(\gamma+2)}$ $(\alpha>0)$, $\s_{6,224}^{\alpha,-\alpha,-1}$ ($\alpha>0$), $\s_{6,227}^{0,-1}$ do not admit complex structures with closed Koszul form can be found in the MAPLE files. 
\end{proof}

\section{Construction of six-dimensional solvmanifolds with holomorphically trivial ca\-no\-ni\-cal bundle via non-invariant sections}\label{S: TCB}
We now work with the Lie algebras that admit complex structures with closed Koszul $1$-form $\psi$, and we construct lattices in the associated simply connected Lie groups such that the canonical bundle of the corresponding solvmanifolds is holomorphically tri\-vial. To achieve this, we first need to give the expression of $\psi$ in each case.

\begin{theorem}\label{theorem: expression of psi}
	Let $\g$ be one of the Lie algebras of Theorem \ref{theorem: CRF}, equipped with a complex structure $J$ such that $\psi([\g,\g])\equiv0$. Then, the expression for $\psi$ in the dual basis $\{e^j\}_{j=1}^6$ of the basis defining the structure equations of the Lie algebra is given in the following table.
	\begin{table}[H]
	\footnotesize{	\begin{tabular}{|c|l|l|} \hline 
			$\g$ & $\psi$ satisfying $\psi([\g,\g])\equiv0$  \\
			\hline 
			$\s_{3,3}^0 \times \R^3$ & $\pm 2 e^3$\\ \hline 
			$\s_{3,3}^0 \times \h_3$ & $\pm 2 e^3$\\ \hline
			$\s_{3,3}^0 \times \s_{3,3}^0$ & $\pm 2 e^3 \pm 2 e^6$ \\ \hline 
			$\s_{4,7} \times \R^2$ & $\pm 2 e^4$\\ \hline 
			$\s_{5,4}^0 \times \R$ & $\pm 2 e^5$\\ \hline 
			$\s_{5,8}^0 \times \R$ & $\pm 4 e^5$ \\ \hline 
			$\s_{5,11}^{\alpha,\alpha,-\alpha}\times \R$ $(\alpha>0)$ & $\pm 2 e^5$ \\ \hline 
			$\s_{5,13}^{\alpha,-\alpha,\gamma} \times \R$ $(0<\gamma\leq1, 0\leq \alpha)$ & $\pm 2 (\gamma\pm 1)e^5$ \\ \hline 
			$\s_{5,16} \times \R$ & $\pm 2 e^5$ \\ \hline 
			$\s_{6,25}$ & $\pm 2 e^6$ \\ \hline 
			$\s_{6,44}$ & $\pm 4 e^6$ \\ \hline 
			$\s_{6,52}^{0,\beta}$ $(\beta>0, \beta\neq 1)$ & $\pm2 (\beta\pm1)e^6$\\ \hline 
			$\s_{6,145}^0$ & $\pm 4 e^6$ \\ \hline 
			$\s_{6,154}^0$ & $\pm4 e^6$ \\ \hline 
			$\s_{6,159}$ & $\pm 2 e^6$ \\ \hline 
			$\s_{6,165}^{\alpha}$ $(\alpha>0)$ & $\pm 4 e^6$ \\ \hline 
			$\s_{6,166}^{\alpha}$ $(0<|\alpha|<1)$ & $\pm2(\alpha\pm1)e^6$ \\ \hline 
			$\s_{6,167}$ & $\pm4 e^6$ \\ \hline 
			$\s_{6,217}^{0,0,-1,-1}$ & $\pm 2 e^5$ \\ \hline 
			$\s_{6,226}^{0,\beta,-1}$ $(0<\beta<1)$ & $\pm2(\beta\pm1) e^5$ \\ \hline 
			$\s_{6,228}^{-\beta,\beta,\gamma,-\gamma}$ $(0\leq\beta\leq\gamma, 0<\gamma)$ & $\pm 2e^5 \pm 2e^6$ \\ \hline
		\end{tabular}}
	\end{table}
\end{theorem}

\begin{proof}
	For each Lie algebra $\g$ in the statement, we set $J$ a generic endomorphism of $\g$, represented as a matrix as in \eqref{eq: J}. We impose the conditions $J^2=-\Id$, $N_J\equiv 0$ and $\psi([\g,\g])\equiv 0$, and we obtain an expression of $\psi$ written in the dual basis $\{e^1,\ldots,e^6\}$ corresponding to the basis $\{e_j\}_{j=1}^6$ that defines the structure equations. 

	For instance, for the Lie algebra 
	$\s_{3,3}^0 \times \R^3=(e^{23},-e^{13},0,0,0,0),$ we obtain that $\psi=a_9 e^1-a_3 e^2+(a_2-a_7)e^3$. Then, for $\psi$ to vanish on the commutator we must have $a_3=a_9=0$. After substitution, looking at the vanishing of the expressions
	\begin{align*}
		B_{33}&=a_{15}^2+a_{16}a_{21}+a_{17}a_{27}+a_{18}a_{33}+1,\\
		N_{141}&=-a_{21}(a_2+a_7), \quad N_{151}=-a_{27}(a_2+a_7), \quad  N_{161}=-a_{33}(a_2+a_7), \\
		N_{142}&=a_{21}(a_1-a_8), \quad N_{152}=a_{27}(a_1-a_8), \quad N_{162}=a_{33}(a_1-a_8), 
	\end{align*} we obtain that $a_7=-a_2$ and $a_8=a_1$. At this point we arrive at $N_{131}=2a_1a_2$ and $  N_{132}=-a_1^2+a_2^2-1$, whose vanishing implies $a_2=\pm 1$, and thus $\psi=\pm 2 e^3$. For the remaining 20 Lie algebras, the computations are explained in full detail in the MAPLE files.
\end{proof}

\medskip

A revised classification of the $6$-dimensional solvable (non-nilpotent) strongly unimodular Lie algebras admitting complex structures with non-zero closed $(3,0)$-forms now follows from the classification of Lie algebras admitting complex structures done in \cite{Sroka} for the decomposable cases $3\times 3$ and $4 \times2$, the almost nilpotent case (\cite{FPAA}, \cite{FPAN1}, \cite{FPAN2}), the nilradical $\R^4$ case (Theorem \ref{theorem: complejas}), the classification of those admitting complex structures with closed Koszul $1$-form (Theorem \ref{theorem: CRF}) and Theorem \ref{theorem: expression of psi}. Hence, the revised version of \cite[Theorem 2.8]{FOU} is the following:	

\begin{theorem}
	Let $\g$ be a solvable (non-nilpotent) strongly unimodular Lie algebra of dimension $6$. Then, $\g$ admits a complex structure with a non-zero closed $(3,0)$-form if and only if it is isomorphic to one in the following list:
	\begin{align*}
		\g_1&\simeq \s_{5,9}^{1,-1,-1}\times \R, \quad \g_2^\alpha=\s_{5,13}^{\alpha,-\alpha,1}\times \R \; (\alpha\geq 0), \quad \g_3\simeq \s_{3,1}^{-1} \times \s_{3,3}^0, \quad \g_4 \simeq \s_{6,52}^{0,1}, \quad \g_5\simeq \s_{6,162}^1,\\
		\g_6&\simeq \s_{6,166}^{-1}, \quad \g_7\simeq \s_{6,166}^1, \quad \g_8\simeq \s_{6,226}^{0,1,-1}, \quad \g_9\simeq \s_{6,152}, \quad \g_{10}=\s_{6,147}^0.
	\end{align*} 
	%
\end{theorem}

Next, we will determine whether the non-zero closed $(3,0)$-form $\tau$ of Proposition \ref{proposition: tau explicit} defined on the associated simply connected Lie group $G$ is invariant under the action of some lattice $\Gamma$ of $G$, which will be constructed using Theorem \ref{theorem: yamada}. In each case we will equip the solvmanifold with the invariant complex structure exhibited in Theorem \ref{theorem: CRF}, since choosing another complex structure only results in a change of sign in the Koszul form (due to Theorem \ref{theorem: expression of psi}), so the lattices that can be used to construct a complex solvmanifold with holomorphically trivial canonical bundle via our method remain the same.

In most cases, we obtain a complex solvmanifold $(\Gamma\backslash G, J)$ biholomorphic to a complex solvmanifold $(\Gamma\backslash \tilde{G}, \tilde{J})$ that has an invariant holomorphic section of its canonical bundle. This is achieved by noticing that $\Gamma$ is also a lattice in $\tilde{G}$, where $\tilde{G}$ is obtained by eliminating the action of the rotations, as a particular case of the so-called \textit{$S$-modification} of a Lie group $G$, introduced by \cite{CFK}. 
 
However, there is one notable exception: the Lie algebra $\s_{6,154}^0$. For this Lie algebra, we construct an associated complex solvmanifold with holomorphically trivial canonical bundle, which we can argue that is not biholomorphic to a complex solvmanifold with an invariant holomorphic section of its canonical bundle. Additionally, we also construct some other interesting examples arising from the Lie algebras $\s_{5,8}^0 \times \R$, $\s_{6,44}$, $\s_{6,145}^0$, $\s_{6,165}^{\alpha}$, $\s_{6,167}$ and $\s_{6,228}^{-\beta,\beta,\gamma,-\gamma}$, where it remains unclear whether they are biholomorphic to a complex solvmanifold admitting invariant holomorphic sections of its canonical bundle. 

\medskip
	
$\bullet$ The Lie algebra $\s_{3,3}^0\times \R^3=(e^{23},-e^{13},0,0,0,0)$ can be expressed as the almost abelian Lie algebra $\R e_3 \ltimes_{\ad} \R^5$, where $\ad e_3|_{\R^5}=\smatriz{0&1\\-1&0}\oplus (0)^{\oplus 3}$. Hence, the corresponding simply connected Lie group is $S_{3,3}^0 \times \R^3:=\R\ltimes_\phi \R^5$, where\footnote{We use $A\oplus B$ to denote the block-diagonal matrix $\left[\begin{smallmatrix}A&\\&B\end{smallmatrix}\right]$. This naturally generalizes to the sum of $n$ square matrices.} \[\phi(t)=\exp(t \ad e_3|_{\R^5})=\matriz{\cos t&\sin t\\ -\sin t &\cos t}\oplus \I_3.\] 
	The closed $(3,0)$-form $\tau$ from \eqref{eq: tau} is given by $\tau=\exp(-it) (e^1+ie^2) \wedge (e^3+ie^4) \wedge (e^5+ie^6)$. We can apply Theorem \ref{theorem: yamada} with $t=2\pi$, since for the rational basis $\mathcal{B}=\{e_1,e_2,e_4,e_5,e_6\}$ of $\R^5$ we have that $[\phi(2\pi)]_{\B}=\I_5\in \SL(5,\Z)$. Thus, $\Gamma=2\pi \Z\ltimes_\phi \Z^5=2\pi \Z \times \Z^5$ is a lattice in $S_{3,3}^0 \times \R^3$ and $\tau$ is invariant under the action of $\Gamma$. Therefore, the complex solvmanifold $(\Gamma\backslash (S_{3,3}^{0}\times \R^3), J)$ has holomorphically trivial canonical bundle. However, this complex solvmanifold is biholomorphic to the complex torus $(\Z^6 \backslash \R^6, J)$. Indeed, $\Gamma$ is a lattice in both $S_{3,3}^0\times \R^3$ and $\R^6$, so that the identity map induces a biholomorphism between the complex solvmanifold $(\Gamma\backslash (S_{3,3}^{0}\times \R^3), J)$ and the complex torus $(\Z^6 \backslash \R^6, J)$.

	\medskip 
	
	$\bullet$ The Lie algebra $\s_{3,3}^0 \times \h_3=(e^{23},-e^{13},0,-e^{56},0,0)$ can be written as $\R e_3 \ltimes_{\ad} (\R^2\times \h_3)$, where $\ad e_3|_{\R^2\times \h_3}=\smatriz{0&1\\-1&0}\oplus (0)^{\oplus 3}$. Hence, the corresponding simply connected Lie group is $S_{3,3}^0\times H_3:= \R \ltimes_\phi (\R^2\times H_3)$, with \[\phi(t)=\matriz{\cos t&\sin t\\ -\sin t &\cos t}\oplus \I_3.\] 
	Here and henceforth, $H_3$ is identified with $\R^3$ equipped with the group law induced by the exponential map.
	 
	The closed $(3,0)$-form $\tau$ from \eqref{eq: tau} is $\tau=\exp(-it) (e^1+ie^2) \wedge (e^3+ie^4) \wedge (e^5+ie^6)$. We apply Theorem \ref{theorem: yamada} with $t=2\pi$, since for the rational basis $\mathcal{B}=\{e_1,e_2,e_4,e_5,e_6\}$ of $\R^2\times H_3$ we have that $[\phi(2\pi)]_{\B}=\I_5\in \SL(5,\Z)$. Thus, $\Gamma=2\pi \Z\ltimes_\phi \Gamma_N$, where $\Gamma_N=\exp^{\R^2\times H_3}\text{span}_\Z\{e_1,e_2,e_4,e_5,e_6\}$, is a lattice in $S_{3,3}^0\times H_3$ and $\tau$ is invariant under the action of $\Gamma$. Therefore, the corresponding complex solvmanifold has trivial canonical bundle. 
	
	We show next that the complex solvmanifold $(\Gamma\backslash (S_{3,3}^{0}\times H_3), J)$ is biholomorphic to a complex nilmanifold $((2\pi \Z\times \Gamma_N) \backslash (\R^3\times H_3), \tilde{J})$. Indeed, the group law can be described as follows: let $\mathbf{x}=(t,x,y,z,u,v), \mathbf{x'}=(t',x',y',z',u',v')$ be elements of $\R^6$. The multiplication in $S_{3,3}^0 \times H_3$ is given by:
	\begin{align*}
	\mathbf{x}\cdot \mathbf{x'}=(t+t',x+x'\cos t+y'\sin t, y-x'\sin t+y'\cos t, z+\tfrac12 (uv'-vu'),u+u',v+v')
	\end{align*}
	and in $\R^3\times H_3$ it is given by:
	\begin{align*}
		\mathbf{x}\cdot \mathbf{x'}=(t+t',x+x',y+y',z+\tfrac12 (uv'-vu'),u+u',v+v').  
	\end{align*}
	Then, let $\{e_j\}_{j=1}^6$ and $\{f_j\}_{j=1}^6$ denote the bases of $\s_{3,3}^0 \times \h_3$ and $\R^3 \times \h_3$, respectively. If we consider $\{e_j\}_{j=1}^6$ and $\{f_j\}_{j=1}^6$ as left-invariant vector fields on $S_{3,3}^0 \times H_3 \equiv \R^3 \times H_3 \equiv \R^6$ in terms of the usual coordinate vector fields on $\R^6$ we have that
	\begin{align*}
		e_1&=\cos t \,\partial_x -\sin t\, \partial_y, \quad e_2=\sin t\, \partial_x+\cos t\, \partial_y, \quad e_3=\partial_t, \quad  e_4=\partial_z, \quad e_5=-\tfrac{v}{2} \partial_z+\partial_u,\\ e_6&=\tfrac{u}{2} \partial_z+\partial_v, \quad f_1=\partial_x, \quad f_2=\partial_y, \quad f_3=\partial_t, \quad f_4=\partial_z, \quad f_5=-\tfrac{v}{2} \partial_z+\partial_u, \quad f_6=\tfrac{u}{2} \partial_z+\partial_v.
	\end{align*}
	From here it can be seen the left-invariant complex structures corresponding to $Je_{2i-1}=e_{2i}$ and $\tilde{J}f_{2i-1}=f_{2i}$, $1\leq i\leq 6$, coincide on $\R^6$, that is, there is a complex structure on $\R^6$ which is invariant by the left actions of both $S_{3,3}^0 \times H_3$ and $\R^3 \times H_3$. We have the following commutative diagram:
	\[\begin{tikzcd}
		(S_{3,3}^0 \times H_3, J) \arrow{r}{F=\Id} \arrow[swap]{d}{\pi_1} & (\R^3 \times H_3,\tilde{J}) \arrow{d}{\pi_2} \\
		(\Gamma\backslash (S_{3,3}^0 \times H_3),J) \arrow{r}{\tilde{F}} & (\Gamma\backslash (\R^3 \times H_3),\tilde{J})
	\end{tikzcd}
	\]
	where the canonical projections $\pi_i$ are local biholomorphisms. Using that the complex structures coincide on $\R^6$ and that $\tilde{F}$ is induced by the identity map, it is easy to see that $\tilde{F}$ is a biholomorphism.

	\medskip
	
	$\bullet$ For the Lie algebra $\s_{3,3}^0 \times \s_{3,3}^0=(e^{23},-e^{13},0,e^{56},-e^{46},0)$,
	 we will change to a new basis, in order to use Proposition \ref{proposition: tau explicit}. Take the basis $\{f_1,\ldots,f_6\}$, where $f_1=e_3+e_6$, $f_2=e_6-e_3$, $f_3=e_1$, $f_4=e_2$, $f_5=e_4$ and $f_6=e_5$. In this basis $\s_{3,3}^0 \times \s_{3,3}^0$ can be written as $(\R f_1 \oplus \R f_2)\ltimes_{\ad} \R^4$, where $\ad f_1|_{\R^4}=\smatriz{0&1\\-1&0}^{\oplus 2}$ and $\ad f_2|_{\R^4}=\smatriz{0&-1\\1&0}\oplus\smatriz{0&1\\-1&0}$. The associated simply connected Lie group is $G=\R^2\ltimes_\phi \R^4$, where  \[\phi(t,0)=\matriz{\cos t&\sin t\\-\sin t&\cos t}^{\oplus 2} \quad \text{and}\quad \phi(0,s)=\matriz{\cos s&-\sin s\\\sin s &\cos s}\oplus \matriz{\cos s&\sin s\\-\sin s &\cos s}.\]
	The complex structure $J$ in the basis $\{f_j\}_{j=1}^6$ is given by 
	\begin{equation}\label{eq: new J} 
	Jf_1=f_2,\quad Jf_3=f_4, \quad \text{and} \quad Jf_5=f_6,
	\end{equation} 
	 and satisfies $\psi=4f^1$. The associated closed $(3,0)$-form $\tau$ defined on $G$ as in \eqref{eq: tau} is given by $\exp(-2it)(f^1+if^2) \wedge (f^3+if^4) \wedge (f^5+if^6)$. To construct a splittable lattice $\Gamma$ that leaves the form $\tau$ invariant we must choose $t\in \{\pi, 2\pi\}$. However, in any case, one can check that the splittable lattices we can get by choosing $s\in \{2\pi, \pi, \frac{\pi}{2},\frac{2\pi}{3},\frac{\pi}{3}\}$ are also lattices in the simply connected Lie group $\tilde{G}$ corresponding to $\g_2^0=\s_{5,13}^{0,0,1}\times \R=(e^{25},-e^{15},e^{45},-e^{35},0,0)$, after changing the basis to $\{\tilde{f_j}\}_{j=1}^6=\{e_5+e_6,e_5-e_6,e_2,e_1,e_3,e_4\}$. This Lie algebra admits a non-vanishing left-invariant closed $(3,0)$-form with respect to the complex structure \eqref{eq: new J}. Moreover, identifying $G\equiv \tilde{G}\equiv \R^6$, in terms of the coordinate vector fields on $\R^6$ we have that
	 \begin{align*}
	 	f_1&=\partial_t, \quad f_2=\partial_x, \quad f_3=\cos(s-t)\partial_y+\sin(s-t)\partial_z, \quad f_4=-\sin(s-t)\partial_y +\cos(s-t)\partial_z \\
	 	f_5&=\cos(t+s)\partial_u-\sin(s+t)\partial_v, \quad 
	 	f_6=\sin(t+s)\partial_u+\cos(t+s)\partial_v,\\ 
	 	\tilde{f_1}&=\partial_t, \quad \tilde{f_2}=\partial_x, \quad \tilde{f_3}=\cos(t+s)\partial_y+\sin(t+s)\partial_z, \quad \tilde{f_4}=-\sin(t+s)\partial_y+\cos(t+s)\partial_z, \\
	 	\tilde{f_5}&=\cos(t+s)\partial_u-\sin(t+s)\partial_v, \quad \tilde{f_6}=\sin(t+s)\partial u+\cos(t+s)\partial_v.
	 \end{align*}
	 Therefore, the left-invariant complex structures coincide on $\R^6$ and thus, as in the previous case, the identity map induces a biholomorphism between the corresponding complex solvmanifolds with holomorphically trivial canonical bundle.
	 
	\medskip

	
	$\bullet$ The Lie algebra $\s_{4,7}\times\R^2=(-e^{23},e^{34},-e^{24},0,0,0)$ can be described as the almost nilpotent Lie algebra $\R e_4 \ltimes_{\ad} (\h_3\times \R^2)$, where $\ad e_4|_{\h_3\times \R^2}=(0)\oplus \smatriz{0&1\\-1&0}\oplus (0)^{\oplus 2}$. Hence, the corresponding simply connected Lie group is $S_{4,7}^0\times \R^2:= \R \ltimes_\phi (H_3\times \R^2)$, where \[\phi(t)=(1)\oplus \matriz{\cos t&\sin t\\ -\sin t &\cos t}\oplus \I_2.\] 
	Therefore, the multiplication is given as follows: let $\mathbf{x}=(t,x,y,z,u,v)$ and $\mathbf{x'}=(t',x',y',z',u',v')$ be elements of $S_{4,7}^0\times \R^2\equiv \R^6$, then
	\begin{align*} 
		\mathbf{x}\cdot\mathbf{x'}=(t+t',&x+x'+\tfrac{y}{2}(-\sin(t) y'+\cos(t) z')-\tfrac{z}{2}(\cos (t) y'+\sin (t) z'),\\
		&\quad  y+\cos (t) y'+\sin (t) z', z-\sin (t) y'+\cos (t) z',u+u',v+v').
		\end{align*} 
	The closed $(3,0)$-form $\tau$, as in \eqref{eq: tau}, is $\tau=\exp(-it) (e^1+ie^4) \wedge (e^2+ie^3) \wedge (e^5+ie^6)$. To construct a lattice $\Gamma$ via Theorem \ref{theorem: yamada} such that $\tau$ is invariant under the action of $\Gamma$ we must choose $t=2\pi$, and a lattice indeed exists since for the rational basis $\mathcal{B}=\{e_1,e_2,e_3,e_5,e_6\}$ of $H_3\times \R^2$ we have that $[\phi(2\pi)]_{\B}=\I_5$. Therefore, $\Gamma=2\pi \Z\ltimes_\phi \Gamma_N$, where $\Gamma_N=\exp^{H_3\times \R^2}\text{span}_\Z\{e_1,e_2,e_3,e_5,e_6\}$, is a lattice in $S_{4,7}\times \R^2$. Considering $\{e_i\}_{i=1}^6$ as left-invariant vector fields on $S_{4,7}\times \R^2\equiv \R^6$, in terms of the coordinate vector fields on $\R^6$ we have that
	\begin{align*}
		e_4&=\partial_t, \quad e_1=\partial_x, \quad e_2=-\tfrac12 (y\sin t+z\cos t)\,\partial_x+\cos t\,\partial_y-\sin t\, \partial_z \\
		 e_3&=\tfrac{1}{2}(y\cos t-z\sin t)\,\partial_x+\sin t\, \partial_y+\cos t\,\partial_z, \quad e_5=\partial_u, \quad e_6=\partial_v. 
		\end{align*} 
			Therefore, the left-invariant complex structure corresponding to $\s_{4,7}\times \R^2$ (from Theorem \ref{theorem: CRF}) on $\R^6$ is given by $J\partial_x=\partial_t$, $J\partial_y=\frac{y}{2}\partial_x+\frac{z}{2}\partial_t+\partial_z$ and $J\partial_u=\partial_v$.
	On another hand, considering the elements $\{f_i\}_{i=1}^6$ of $\R f_4 \times \h_3 \times \R^2$ as left-invariant vector fields on $\R^6$, in terms of the coordinate vector fields we have that
	\begin{align*}
		f_4=\partial_t,& \quad f_1=\partial_x, \quad f_2=-\tfrac{z}{2} \partial_x+\partial_y, \quad f_3=\tfrac{y}{2}\partial_x+\partial_z, \quad f_5=\partial_u, \quad f_6=\partial_v.  
	\end{align*}
A straightforward calculation shows that both complex structures coincide on $\R^6$, and since $\Gamma$ is also a lattice in $\R\times H_3\times \R^2$ the identity map $S_{4,7}\times\R^2\to \R\times H_3\times \R^2$ induces a biholomorphism between  $(\Gamma\backslash (S_{4,7}\times \R^2), J)$ and the complex nilmanifold $((2\pi \Z\times \Gamma_N) \backslash (\R\times H_3 \times \R^2), J)$. 

\medskip
		
	$\bullet$ The Lie algebra $\s_{5,4}^0 \times \R=(e^{25},0,e^{45},-e^{35},0,0)$ can be written
	as $\R e_5\ltimes_{\ad} \R^5$, where $\ad e_5|_{\R^5}=\smatriz{0&1\\0&0}\oplus\smatriz{0&1\\-1&0}\oplus (0)$, and hence the corresponding simply connected Lie group is $S_{5,4}^0\times \R:= \R \ltimes_\phi \R^5$, where 
	\[  \phi(t)=\matriz{1&t\\0&1}\oplus \matriz{\cos t&\sin t\\ -\sin t&\cos t}\oplus (1).\] 
	The closed $(3,0)$-form $\tau$ from \eqref{eq: tau} is $\tau=\exp(-it)(e^1+ie^6) \wedge (e^2+ie^5) \wedge (e^3+ie^4)$. We apply Theorem \ref{theorem: yamada} with $t=2\pi$, since for the rational basis $\B=\{2\pi e_1, e_2,e_3,e_4,e_6\}$ of $\R^5$ we have that $[\phi(2\pi)]_\B=\smatriz{1&1\\0&1}\oplus \I_3$. Consequently, $\Gamma=2\pi \Z\ltimes_\phi  \text{span}_\Z\{2\pi e_1, e_2, e_3, e_4, e_6\}$ is a lattice in $S_{5,4}^0 \times \R$ and $\tau$ is invariant under the action of $\Gamma$. However, as in previous cases, one can check that the complex solvmanifold with holomorphically trivial canonical bundle $(\Gamma\backslash (S_{5,4}^0 \times \R),J)$ is biholomorphic to the complex nilmanifold $(\Gamma \backslash N, J)$, where $N$ is the nilpotent, simply connected Lie group associated to the Lie algebra $\n:=\text{span}\{e_1,\ldots,e_6\}$, with the only non-trivial Lie bracket being $[e_5,e_2]=e_1$. 
	
	\medskip
	
	$\bullet$ The Lie algebra $\s_{5,8}^0 \times \R=(e^{25}+e^{35},-e^{15}+e^{45},e^{45},-e^{35},0,0)$ can be written as the almost abelian Lie algebra $\R e_5 \ltimes_{\ad} \R^5$, where $\ad e_5|_{\R^5}=\left[\begin{smallmatrix}0&1&1&0\\-1&0&0&1\\0&0&0&1\\0&0&-1&0\end{smallmatrix}\right]\oplus (0)$. Hence, the associated simply connected Lie group is $S_{5,8}^0 \times \R:=\R\ltimes_\phi \R^5$, where \[\phi(t)=\matriz{\cos t&\sin t& t \cos t & t \sin t\\ -\sin t & \cos t & -t \sin t& t\cos t\\0&0&\cos t&\sin t\\0&0&-\sin t&\cos t}\oplus (1).\] The closed non-zero $(3,0)$-form $\tau$ from \eqref{eq: tau} is $\tau=\exp(-2it)(e^1+ie^2)\wedge (e^3+ie^4)\wedge (e^5+ie^6)$. In order to construct a lattice $\Gamma$ via Theorem \ref{theorem: yamada} such that $\tau$ is invariant under the action of $\Gamma$ we must choose $t\in \{2\pi, \pi\}$. For $t=2\pi$, $\B=\{2\pi e_1, e_3, 2\pi e_2, e_4, e_6\}$ is a rational basis of $\R^5$ satisfying $[\phi(2\pi)]_\B=\left[\begin{smallmatrix}1&1\\0&1\end{smallmatrix}\right]^{\oplus 2}\oplus 1$. According to Theorem \ref{theorem: yamada}, $\Gamma=2\pi \Z\ltimes_\phi \text{span}_\Z \B$ is a lattice in $S_{5,8}^0\times \R$. However, as in the previous cases, one can check that the associated complex solvmanifold with holomorphically trivial canonical bundle is biholomorphic to one associated to the nilpotent Lie algebra $\n=\text{span}\{e_1,\ldots,e_6\}$ with non-trivial Lie brackets given by $[e_5, e_3]=e_1, [e_5, e_4]=e_2$, equipped with the same complex structure and quotienting by the same lattice.
	
	For $t=\pi$ we have that $[\phi(\pi)]_\B=\left[\begin{smallmatrix}-1&1\\0&-1\end{smallmatrix}\right]^{\oplus 2}\oplus (1)$, where $\B$ is the rational basis of $\R^5$ given by $\B=\{\pi e_1, -e_3, \pi e_2, -e_4, e_6\}$. Hence, $\Gamma=\pi \Z\ltimes_{\phi} \text{span}_\Z\{\pi e_1, \pi e_3, -e_2, e_4, e_6\}$ is a lattice in $S_{5,8}^0 \times \R$. The corresponding complex solvmanifold $(\Gamma\backslash S_{5,8}^0,J)$ has holomorphically trivial canonical bundle.
	
	
	\medskip
	
	$\bullet$ The Lie algebra $\s_{5,11}^{\alpha, \alpha, -\alpha}\times \R=(\alpha e^{15}, \alpha e^{25}, -\alpha  e^{35}+e^{45}, -e^{35}-\alpha e^{45},0,0)$, where $\alpha>0,$ can be expressed as $\R e_5 \ltimes_{\ad} \R^5$, where $\ad e_5|_{\R^5}=(\alpha)^{\oplus 2}\oplus \left[\begin{smallmatrix}-\alpha&1\\-1&-\alpha\end{smallmatrix}\right]\oplus (0)$. Hence, the corresponding simply connected Lie group is $S_{5,11}^{\alpha,\alpha,-\alpha}\times \R:=\R\ltimes_\phi \R^5$, where \[ \phi(t)=(e^{\alpha t})^{\oplus 2}\oplus e^{-\alpha t} \matriz{\cos t&\sin t\\-\sin t&\cos t}\oplus (1).\] 
The closed $(3,0)$-form $\tau$ from \eqref{eq: tau} is $\tau=\exp(-it)(e^1+ie^2) \wedge (e^3+ie^4) \wedge (e^5+ie^6)$. To construct a lattice $\Gamma$ via Theorem \ref{theorem: yamada} such that $\tau$ is invariant under the action of $\Gamma$ we must choose $t=2\pi$. If we set $\alpha=\frac{t_m}{2\pi}$, where $t_m=\log(\frac{m+\sqrt{m^2-4}}{2})$, for $m\in \N$, $m\geq 3$, then the basis $\B=\{e_1+e_3, \exp(t_m)e_1+\exp(-t_m)e_3, e_2+e_4, \exp(t_m)e_2+\exp(-t_m)e_4,e_6\}$ is a rational basis of $\R^5$ and satisfies  $[\phi(2\pi)]_\B=\left[\begin{smallmatrix}0&-1\\1&m\end{smallmatrix}\right]^{\oplus 2}\oplus (1)\in \SL(5,\Z)$. Hence, $\Gamma_m=2\pi \Z\ltimes_\phi \text{span}_\Z \mathcal{B}$ is a lattice in $S_{5,11}^{\alpha,-\alpha,-\alpha}\times \R$. However, it can be seen that the associated complex solvmanifold is biholomorphic to the complex solvmanifold $(\Gamma_m\backslash (S_{5,9}^{1,-1,-1}\times \R), J)$, using that $\Gamma_m$ is also a lattice in the simply connected Lie group $S_{5,9}^{1,-1,-1}\times \R$ and that the corresponding left-invariant complex structures coincide on $\R^6$. Here $S_{5,9}^{1,-1,-1}\times \R$ is the simply connected Lie group associated to $\s_{5,9}^{1,-1,-1}\times \R$, which admits a non-zero closed $(3,0)$-form with respect to the complex structure $J$.

	\medskip
	
	$\bullet$ The Lie algebra $\s_{5,13}^{\alpha, -\alpha, \gamma}\times \R=(\alpha e^{15}+e^{25},-e^{15}+\alpha e^{25}, -\alpha e^{35}+\gamma e^{45}, -\gamma e^{35}-\alpha e^{45},0,0)$, where $0<\gamma<1$ and $0\leq \alpha$, can be expressed as the almost abelian Lie algebra $\R e_5\ltimes_{\ad e_5} \R^5$, where $\ad e_5|_{\R^5}=\smatriz{\alpha&1\\-1&\alpha}\oplus\smatriz{-\alpha&\gamma\\-\gamma&-\alpha}\oplus (0)$, and the corresponding simply connected Lie group is $S_{5,13}^{\alpha,-\alpha,\gamma}\times \R:=\R\ltimes_\phi \R^5$, where 
	\[\phi(t)=e^{\alpha t}\matriz{\cos t&\sin t\\ -\sin t & \cos t} \oplus e^{-\alpha t}\matriz{\cos(\gamma t)&\sin(\gamma t)\\ -\sin(\gamma t)&\cos(\gamma t)}\oplus (1).\]
	The closed $(3,0)$-form $\tau$ from \eqref{eq: tau} is $\tau=\exp(-i(\gamma-1)t)(e^1-ie^2)\wedge (e^3+ie^4)\wedge (e^5+ie^6)$. The only way to obtain a lattice $\Gamma$ using Theorem \ref{theorem: yamada} such that the form $\tau$ is invariant under the action of $\Gamma$ is to choose $t=\frac{2\pi k}{\gamma-1}$ for some $k\in \Z\setminus\{0\}$, but this forces $\gamma t\in t+2\pi \Z$. Hence, the corresponding lattices will be also lattices in the solvable group $S_{5,13}^{\alpha,-\alpha,1}\times \R$, which admits a non-vanishing left-invariant closed $(3,0)$-form with respect to the same complex structure. Thus, it can be seen as in previous cases that the identity map produces a biholomorphism between the complex solvmanifolds. 
	
	Note that a lattice does not exist for every choice of $\alpha$ and $\gamma$. Indeed, according to a result of \cite{Witte}, in each dimension there are only countably many Lie groups that admit lattices. A possible choice to get a lattice using Theorem \ref{theorem: yamada} is to take $\gamma=\frac{1}{n}$ with $n\in \N$ and $\alpha=\frac{t_m}{2\pi n}$, with $t_m$ defined as we did for $S_{5,11}^{\alpha,\alpha,-\alpha}\times \R$. For $t=2\pi n$ we obtain $\phi(t)=e^{t_m} \Id_2 \oplus e^{-t_m} \Id_2 \oplus (1)$, which can be conjugated to an integer matrix as above.

	\medskip
	
	$\bullet$ The Lie algebra $\s_{5,16}\times \R=(-e^{23}+e^{45},e^{35},-e^{25},0,0,0)$ can be expressed as $\R e_5\ltimes_{\ad} (\h_3 \times \R^2)$, where $\ad e_5|_{\h_3\times \R^2}=\left[\begin{smallmatrix} 0&0&0&1\\0&0&1&0\\0&-1&0&0\\0&0&0&0 \end{smallmatrix}\right]\oplus (0)$. Hence, the corresponding simply connected Lie group is $S_{5,16}\times \R:=\R \ltimes_\phi (H_3\times \R^2)$, where \[\phi(t)=\matriz{1&0&0&t\\0&\cos t&\sin t&0\\0&-\sin t&\cos t&0\\0&0&0&1}\oplus (1).\] The closed $(3,0)$-form as in \eqref{eq: tau} is $\tau=\exp(-it)(e^1+ie^6)\wedge (e^2+ie^3) \wedge (e^4+ie^5)$. Hence, as in previous cases, we are forced to choose $t=2\pi$. Let $\B$ be the rational basis of $\h_3\times \R^2$ given by $\B=\{2\pi e_1, e_4, 2\pi e_2, e_3, e_6\}$. Then, given that $[\phi(2\pi)]_\B=\smatriz{1&1\\0&1}\oplus \I_3\in \SL(5,\Z)$, it follows from Theorem \ref{theorem: yamada} that $\Gamma=2\pi \Z \ltimes_\phi \Gamma_N$, where $\Gamma_N=\exp^{H_3\times \R^2}\text{span}_\Z \B$, is a lattice in $S_{5,16}\times \R$. As in previous cases, the complex solvmanifold with holomorphically trivial canonical bundle $(\Gamma\backslash (S_{5,16}\times \R),J)$ is biholomorphic to the complex nilmanifold $(\Gamma\backslash N,J)$, where $N$ is the nilpotent simply connected Lie group corresponding to the Lie algebra spanned by $\{e_1,\ldots, e_6\}$ such that $[e_2,e_3]=e_1, [e_4,e_5]=-e_1$.
	
	\medskip
	
	$\bullet$ The Lie algebra $\s_{6,25}=(-e^{23},e^{36},-e^{26},0,e^{46},0)$ can be written as $\R e_6 \ltimes_{\ad e_6} (\h_3 \times \R^2)$, where $\ad e_6|_{\h_3\times \R^2}=(0)\oplus \smatriz{0&1\\-1&0}\oplus\smatriz{0&0\\1&0}$, and thus the corresponding simply connected Lie group is $S_{6,25}:=\R\ltimes_\phi(H_3\times \R^2)$, where \[ \phi(t)=(1)\oplus \matriz{\cos t&\sin t\\ -\sin t& \cos t}\oplus \matriz{1&0\\t&1}.\] The closed non-zero $(3,0)$-form from \eqref{eq: tau} is $\tau=\exp(-it)(e^1+ie^5)\wedge (e^2+ie^3) \wedge (e^4+ie^6)$. As in previous cases, we must choose $t=2\pi$. Since the rational basis $\B=\{e_1,e_2,e_3,2\pi e_5,e_4+e_5\}$ of $H_3\times \R^2$ satisfies that $[\phi(2\pi)]_\B=\I_3\oplus\smatriz{1&1\\0&1}\in \SL(5,\Z)$, we can apply Theorem \ref{theorem: yamada} to get that $\Gamma=2\pi \Z\ltimes_\phi \exp^{H_3\times \R^2} \text{span}_\Z \B$ is a lattice in $S_{6,25}$ that leaves $\tau$ invariant. In a completely analogous way to the previous cases, one can check that the complex solvmanifold with holomorphically trivial canonical bundle $(\Gamma\backslash S_{6,25},J)$ is biholomorphic to the complex nilmanifold $(\Gamma\backslash N,J)$, where $N$ is the nilpotent simply connected Lie group corresponding to the Lie algebra spanned by $\{e_j\}_{j=1}^6$ with non-trivial Lie brackets $[e_2,e_3]=e_1, [e_4,e_6]=-e_5$. 		
	
	\medskip
	
$\bullet$ The Lie algebra $\s_{6,44}=(-e^{23},e^{36},-e^{26},e^{26}+e^{56},e^{36}-e^{46},0)$ can be written as the almost nilpotent Lie algebra $\R e_6\ltimes_{\ad} (\h_3\times \R^2)$ with $\ad e_6|_{\h_3\times \R^2}=(0)\oplus \left[\begin{smallmatrix}
	0&1&0&0\\ -1&0 &0&0\\ 1&0&0&1\\ 0&1&-1&0
\end{smallmatrix}\right]$, and thus the corresponding simply connected Lie group is $S_{6,44}:=\R\ltimes_\phi (H_3\times \R^2)$, where \[\phi(t)=(1)\oplus \matriz{\cos t&\sin t&0&0\\ -\sin t& \cos t&0&0\\ t \cos t&t \sin t&\cos t&\sin t\\ -t \sin t& t \cos t&-\sin t& \cos t}.\]
The closed $(3,0)$-form $\tau$ from \eqref{eq: tau} is $\tau=\exp(-2it)(e^1+ie^6)\wedge (e^2+ie^3)\wedge (e^4+ie^5)$. To construct a lattice $\Gamma$ via Theorem \ref{theorem: yamada} such that $\tau$ is invariant under the action of $\Gamma$ we are forced to choose $t\in \{\pi, 2\pi\}$. For $t=2\pi$ one can check that we obtain a complex solvmanifold with holomorphically trivial canonical bundle biholomorphic to a complex nilmanifold arising from the nilpotent Lie algebra spanned by $\{e_j\}_{j=1}^6$ whose only non-trivial Lie brackets are $[e_2,e_3]=e_1, [e_3,e_6]=-e_5, [e_2,e_6]=-e_4$.

For $t=\pi$, we have $[\phi(\pi)]_\B=(1)\oplus\smatriz{-1&1\\0&-1}^{\oplus 2}$, where $\B=\{ e_1, -\pi e_4, e_2, -\pi e_5, e_3\}$ is a rational basis of $\h_3\times \R^2$. According to Theorem \ref{theorem: yamada},  $\Gamma=\pi \Z\ltimes_\phi \exp^{H_3\times \R^2} \text{span}_\Z \B$ is a lattice in $S_{6,44}$. The corresponding complex solvmanifold has holomorphically trivial canonical bundle.

	\medskip
	
$\bullet$ The Lie algebra $\s_{6,52}^{0,\beta}=(-e^{23},e^{36},-e^{26},\beta e^{56},-\beta e^{46},0)$, where $\beta>0, \beta\neq 1$, can be expressed as $\R e_6 \ltimes_{\ad} (\h_3\times \R^2)$, where $\ad e_6|_{\h_3\times \R^2}=(0)\oplus\smatriz{0&1\\-1&0}\oplus \smatriz{0&\beta\\-\beta&0}$. Hence, the associated simply connected Lie group is $S_{6,52}^{0,\beta}:=\R\ltimes_\phi (H_3\times \R^2)$ where \[ \phi(t)=(1)\oplus \matriz{\cos t&\sin t\\ -\sin t&\cos t}\oplus \matriz{\cos(\beta t)&\sin(\beta t)\\ -\sin(\beta t)&\cos(\beta t)}.\] The closed $(3,0)$-form $\tau$ from \eqref{eq: tau} is $\tau=\exp(-i(\beta+1)t)(e^1+ie^6)\wedge (e^2+ie^3) \wedge (e^4+ie^5)$. To obtain a lattice $\Gamma$ using Theorem \ref{theorem: yamada} such that the form $\tau$ is invariant under the action of $\Gamma$, our only option is to choose $t=\frac{2\pi k}{\beta+1}$ for some $k\in \Z\setminus\{0\}$, ensuring that $\beta t\in -t+2\pi \Z$. Consequently, the possible lattices we can obtain by choosing $\beta$ and $k$ such that $t\in\{2\pi, \pi, \frac{\pi}{2}, \frac{2\pi}{3}, \frac{\pi}{3}\}+2\pi\Z$ can also be obtained by considering $S_{6,52}^{0,1}$, which admits a non-vanishing left-invariant closed $(3,0)$-form. Again it can be seen that the corresponding complex solvmanifolds with holomorphically trivial canonical bundle are biholomorphic to complex solvmanifolds admitting an invariant trivializing holomorphic section of their canonical bundle. 

\medskip
	
	$\bullet$ The Lie algebra $\s_{6,145}^0=(e^{26}-e^{35},-e^{16}-e^{45},e^{46},-e^{36},0,0)$ can be expressed as $\R e_6\ltimes_{\ad} \n_{5,1}$, with $\ad e_6|_{\n_{5,1}}=\smatriz{0&1\\-1&0}^{\oplus 2}\oplus (0)$. The associated simply connected Lie group is $S_{6,145}^0:=\R\ltimes_\phi N_{5,1}$, where 
	\[ \phi(t)=\matriz{\cos t&\sin t\\ -\sin t&\cos t}^{\oplus 2}\oplus (1).\] The closed non-zero $(3,0)$-form from \eqref{eq: tau} is $\tau=\exp(-2it)(e^1+ie^2)\wedge (e^3+ie^4)\wedge (e^5+ie^6)$. To construct a lattice $\Gamma$ using Theorem \ref{theorem: yamada} such that $\tau$ is $\Gamma$-invariant, we must choose $t\in\{2\pi, \pi\}$. As in previous cases, for $t=2\pi$, the associated complex solvmanifold with holomorphically trivial canonical bundle will be biholomorphic to a complex nilmanifold arising from the Lie algebra $\R\times \n_{5,1}$.
	
	If we choose $t=\pi$, then we have that the rational basis $\B=\{e_1,\ldots,e_5\}$ of $\n_{5,1}$ satisfies  $[\phi(\pi)]_\B=-\I_4\oplus (1)\in \SL(5,\Z)$. Therefore $\Gamma=\pi \Z\ltimes_\phi \exp^{N_{5,1}} \text{span}_\Z \B$ is a lattice in $S_{6,145}^0$ and the corresponding complex solvmanifold has holomorphically trivial canonical bundle. 
	
	
	\medskip
	
	$\bullet$ The Lie algebra $\s_{6,154}^0=(-e^{26}-e^{35},e^{16}-e^{34},-e^{45},e^{56},-e^{46},0)$ can be written as $\R e_6\ltimes_{\ad} \n_{5,2}$, where $\ad e_6|_{\n_{5,2}}=\smatriz{0&-1\\1&0}\oplus (0)\oplus \smatriz{0&1\\-1&0}$. Hence, the associated simply connected Lie group is $S_{6,154}^0:=\R\ltimes_\phi N_{5,2}$, where \[\phi(t)=\matriz{\cos t&-\sin t\\ \sin t&\cos t} \oplus (1)\oplus \matriz{\cos t&\sin t\\ -\sin t&\cos t}.\] The closed $(3,0)$-form $\tau$ from \eqref{eq: tau} is $\tau=\exp(-2it)(e^1-ie^2)\wedge (e^3+ie^6)\wedge (e^4+ie^5)$. As in previous cases, we are forced to choose $t\in\{2\pi, \pi\}$. For $t=2\pi$, the discrete subgroup $\Gamma_{2\pi}:=2\pi \Z\ltimes_\phi \exp^{N_{5,2}}\text{span}_\Z\{e_1,\ldots,e_5\}$ of $S_{6,154}^0$ is a lattice and, since $\Gamma_{2\pi}$ is also a lattice in $\R \times N_{5,2}$, the complex solvmanifold $M_{2\pi}:=\Gamma_{2\pi}\backslash S_{6,154}^0$ with holomorphically trivial canonical bundle is biholomorphic to the complex nilmanifold $\Gamma_{2\pi} \backslash (\R\times N_{5,2})$.
	
	If we choose $t=\pi$, then we have that the rational basis $\B=\{e_1,\ldots,e_5\}$ of $\n_{5,2}$ satisfies  $[\phi(\pi)]_\B=-\I_2\oplus (1)\oplus -\I_2$ and hence $\Gamma_{\pi}:=\pi \Z\ltimes_\phi \exp^{N_{5,2}} \text{span}_\Z \{e_1,\ldots,e_5\}$ is a lattice in $S_{6,154}^0$. 
	
	We compute next the first Betti number of the corresponding complex solvmanifold with holomorphically trivial canonical bundle, since it will be used later (see Proposition \ref{prop: S6,154} below).
	We identify $\Gamma_\pi$ with $\Z\ltimes_{[\phi(\pi)]_\B} \Gamma_N$, where $\Gamma_N=\text{span}_\Z \{e_1,\ldots,e_5\}$ with the product given by the BCH formula: $x\cdot y=x+y+\tfrac12 [x,y]+\tfrac{1}{12}([x,[x,y]]+[y,[y,x]])$. Thus, the product of $\Gamma_N$ is given by
	\begin{align*} 
		\left(\sum_{i=1}^5 m_i e_i\right) \cdot \left(\sum_{i=1}^5 n_i e_i\right)&=\sum_{i=1}^5 (m_i+n_i)e_i+(\tfrac{m_4 n_5^2-m_5 n_4n_5-m_4m_5n_5+m_5^2 n_4}{12}+\tfrac{m_3n_5-m_5n_3}{2})e_1\\
		&+\quad (\tfrac{m_4n_4n_5-m_5n_4^2-m_4^2n_5+m_4m_5n_4}{12}+\tfrac{m_3n_4-m_4n_3}{2})e_2+(\tfrac{m_4n_5-m_5n_4}{2})e_3 .
		\end{align*}
	 A set of generators of $\Gamma_N$ is $\{\frac{e_1}{6}, \frac{e_2}{6}, \frac{e_3}{2}, e_4,e_5\}$, and the only non-trivial relations are $e_3^{e_4}=e_2\cdot e_3$, $e_3^{e_5}=e_1\cdot e_3$ and $e_4^{e_5}=(\frac{e_1}{2})\cdot (-\frac{e_2}{2})\cdot e_3\cdot e_4$. Therefore, 
	\begin{align*}  
		\Gamma_\pi\cong \la t_1,\ldots, t_6 \mid& t_3^{t_4}=t_2^3 t_3, \; t_3^{t_5}=t_1^3 t_3,\; t_4^{t_5}=t_1^3 t_2^{-3} t_3^2 t_4,\\ & t_1^{t_6}=t_1^{-1}, \, t_2^{t_6}=t_2^{-1},\, t_3^{t_6}=t_3,\, t_4^{t_6}=t_4^{-1}, t_5^{t_6}=t_5^{-1}\ra.
		\end{align*}
	Here $t_i^{t_j}:=t_j^{-1} t_i t_j$. Hence, the abelianization of $\Gamma_\pi$ is $\Z_2^3\oplus \Z$, so that $b_1(\Gamma\backslash S_{6,154}^0)=1$.  
	
	\begin{remark}\label{rmk: 154 no symplectic}
		We observe that $M_{\pi}:=\Gamma_{\pi}\backslash S_{6,154}^0$ does not admit any symplectic form. Assume, for the sake of contradiction, that $M_{\pi}$ admits a symplectic form $\omega$ and consider the map $p:M_{2\pi}\to M_{\pi}$ defined by $p([g]_{\Gamma_{2\pi}})=[g]_{\Gamma_{\pi}}$. This map is well-defined since $\Gamma_{\pi}$ is a subgroup of $\Gamma_{2\pi}$, and it is a surjective submersion, as it satisfies $p\circ p_{2\pi}=p_{\pi}$, where $p_{2\pi}$ and $p_{\pi}$ denote the natural projections $S_{6,154}^0\to M_{2\pi}$ and $S_{6,154}^0\to M_{\pi}$ respectively. 
		
		Since $p$ is a surjective submersion, the pullback $p^* \omega$ would be a symplectic form on $M_{2\pi}$, because $d(p^*\omega)=p^*(d\omega)=0$ and $(p^* \omega)^3=p^* \omega^3\neq 0$, as $p^*$ is injective and $\omega^3\neq 0$.
		
		Now, since $M_{2\pi}$ is diffeomorphic to the nilmanifold $N_{2\pi}:=\Gamma_{2\pi} \backslash (\R\times N_{5,2})$, this would imply that $N_{2\pi}$ admits a symplectic form. By Nomizu's theorem, the de Rham cohomology of $N_{2\pi}$ can be computed using invariant differential forms, so the corresponding Lie algebra $\R\times \n_{5,2}$, would also admit a symplectic form. However, a straightforward calculation shows that this Lie algebra does not admit any symplectic form, leading to a contradiction. Therefore, we can conclude that  $M_{\pi}=\Gamma_\pi \backslash S_{6,154}^0$ does not admit any symplectic form.
	\end{remark}

	$\bullet$ The Lie algebra $\s_{6,159}=(-e^{24}-e^{35},0,-e^{56},0,e^{36},0)$ can be expressed as $\R e_6\ltimes_{\ad e_6} \n_{5,3}$, where $\ad e_6|_{\n_{5,3}}=(0)^{\oplus 2}\oplus \left[\begin{smallmatrix}0&0&-1\\0&0&0\\1&0&0\end{smallmatrix}\right]$. Thus, the corresponding simply connected Lie group is $S_{6,159}:=\R\ltimes_\phi N_{5,3}$, where \[ \phi(t)=\I_2\oplus \matriz{\cos t&0&-\sin t\\0&1&0\\\sin t&0&\cos t}.\] The closed $(3,0)$-form from \eqref{eq: tau} is given by $\tau=\exp(it)(e^1+ie^6)\wedge (e^2+ie^4) \wedge (e^3+ie^5)$. Hence, we must choose $t=2\pi$, but then the associated complex solvmanifold with holomorphically trivial canonical bundle will be biholomorphic to a complex nilmanifold arising from $\R\times \n_{5,3}$.  
	
	\medskip
	
	
	$\bullet$ The Lie algebra $\s_{6,165}^{\alpha}=(-e^{24}-e^{35},\alpha e^{26}-e^{36},e^{26}+\alpha e^{36}, -\alpha e^{46}-e^{56}, e^{46}-\alpha e^{56},0)$, where $\alpha>0$, can be written as $\R e_6\ltimes_{\ad} \n_{5,3}$ with $\ad e_6|_{\n_{5,3}}=(0)\oplus \smatriz{\alpha&-1\\1&\alpha}\oplus \smatriz{-\alpha&-1\\1&-1}$. Thus, the associated simply connected Lie group is $S_{6,165}^\alpha:=\R\ltimes_\phi N_{5,3}$, where \[ \phi(t)=(1)\oplus e^{\alpha t} \matriz{\cos t&-\sin t\\ \sin t& \cos t}\oplus e^{-\alpha t} \matriz{\cos t & -\sin t\\ \sin t & \cos t}.\]
The closed $(3,0)$-form $\tau$ from \eqref{eq: tau} is $\tau=\exp(2it)(e^1+ie^6)\wedge (e^2+ie^3)\wedge (e^4+ie^5)$. To obtain a lattice $\Gamma$ using Theorem \ref{theorem: yamada} such that the form $\tau$ is $\Gamma$-invariant we must choose $t\in \{2\pi, \pi\}$. Moreover, there are only countably many values of $\alpha$ for which a lattice exists. By choosing $t=2\pi$ (and $\alpha=\frac{t_m}{2\pi}$ as below), the resulting complex solvmanifold with holomorphically trivial canonical bundle is biholomorphic to a complex solvmanifold associated with $\s_{6,162}^1$, which admits a non-zero closed $(3,0)$-form. 

For $t=\pi$, we set $\alpha=\frac{t_m}{\pi}$,  where $t_m=\log\left(\frac{m+\sqrt{m^2-4}}{2}\right), m\geq 3$. We have that the rational basis $\B_m=\{f_1,\ldots,f_5\}=\{(\operatorname{e}^{t_m}-\operatorname{e}^{-t_m})e_1, e_2+e_4, -\operatorname{e}^{t_m}e_2-\operatorname{e}^{-t_m}e_4, e_3+e_5, -\operatorname{e}^{t_m}e_3-\operatorname{e}^{-t_m}e_5\}$ of $\n_{5,3}$ satisfies that the only non-trivial Lie brackets are $[f_2,f_3]=f_1$ and $[f_4,f_5]=f_1$, and moreover, $[\phi(\pi)]_\B=(1)\oplus \left[\begin{smallmatrix}0&-1\\1&-m \end{smallmatrix}\right]^{\oplus 2}$. Thus, $\Gamma_m=\pi \Z\ltimes_\phi \exp^{N_{5,3}}\text{span}_\Z\B_m$ is a lattice in $S_{6,165}^0$. The corresponding complex solvmanifold has holomorphically trivial canonical bundle.


\medskip

	$\bullet$ The Lie algebra $\s_{6,166}^{\alpha}=(-e^{24}-e^{35},-e^{46},-\alpha e^{56},e^{26}, \alpha e^{36},0)$, where $0<|\alpha|<1$, can be written as $\R e_6\ltimes_{\ad} \n_{5,3}$, with $\ad e_6|_{\n_{5,3}}=(0)\oplus \left[\begin{smallmatrix}
		0&0&-1&0\\
		0&0&0&-\alpha\\
		1&0&0&0\\
		0&\alpha&0&0
	\end{smallmatrix}\right]$. The associated simply connected Lie group is $S_{6,166}^\alpha:=\R\ltimes_\phi N_{5,3}$, where 
	\[ \phi(t)=(1)\oplus \matriz{\cos t &0 &-\sin t&0 \\ 0 &\cos(t\alpha)&0&-\sin(t\alpha)\\ \sin t &0 & \cos t &0\\ 0 &\sin(t\alpha)&0&\cos(t\alpha)}  .\] 
	The closed $(3,0)$-form $\tau$ from \eqref{eq: tau} is $\exp(i(\alpha+1)t)(e^1+ie^6)\wedge (e^2+ie^4) \wedge (e^3+ie^5)$. To construct a lattice $\Gamma$ using Theorem \ref{theorem: yamada} (for the countably many values of $\alpha$ for which a lattice exists) such that $\tau$ is $\Gamma$-invariant we are forced to choose $t=\frac{2\pi k}{\alpha+1}$ with $k\in \Z\setminus\{0\}$, ensuring that $t\alpha\in -t+2\pi \Z$. Therefore, the corresponding lattices can also be obtained by considering the simply connected Lie group $\s_{6,166}^{-1}$, which admits a non-zero closed $(3,0)$-form with respect to the same complex structure we equipped $\s_{6,166}^{-1}$ in Theorem \ref{theorem: CRF}, and as in previous cases the associated complex solvmanifolds with holomorphically trivial canonical bundle will be biholomorphic. 
	
	To get lattices using Theorem \ref{theorem: yamada}, we can choose, for instance, $\alpha=\frac{1}{n}$ with $n\in \N, n\geq 2$ and $t=2\pi n$, since in the rational basis $\mathcal{B}=\{e_1,e_2,e_3,e_4,e_5\}$ of $\n_{5,3}$ we get $[\phi(2\pi n)]_{\mathcal{B}}=\Id_5\in \SL(5,\Z)$.
	\medskip
		
	$\bullet$ The Lie algebra $\s_{6,167}=(-e^{24}-e^{35},-e^{36},e^{26},e^{26}-e^{56},e^{36}+e^{46},0)$ can be written as the almost nilpotent Lie algebra $\R e_6\ltimes_{\ad} \n_{5,3}$ with $\ad e_6|_{\n_{5,3}}=(0)\oplus \left[ \begin{smallmatrix} 0&-1&0&0\\1&0&0&0\\1&0&0&-1\\0&1&1&0 \end{smallmatrix}\right]$. Thus, the associated simply connected Lie group is $S_{6,167}:=\R\ltimes_\phi N_{5,3}$, where
	\[ \phi(t)=(1)\oplus \matriz{\cos t&-\sin t&0&0\\ \sin t & \cos t &0&0\\ t \cos t& -t \sin t&\cos t&-\sin t\\ t \sin t& t \cos t&\sin t&\cos t}.\]
	The closed $(3,0)$-form $\tau$ from \eqref{eq: tau} is $\tau=\exp(2it)(e^1+ie^6)\wedge(e^2+ie^3)\wedge(e^4+ie^5)$. In order to obtain a lattice $\Gamma$ using Theorem \ref{theorem: yamada} such that $\tau$ is $\Gamma$-invariant, we are forced to choose $t\in \{2\pi, \pi\}$. Similarly to previous cases, for $t=2\pi$ the corresponding complex solvmanifold with holomorphically trivial canonical bundle will be biholomorphic to a complex nilmanifold arising from the nilpotent Lie algebra spanned by $\{e_1,\ldots, e_6\}$ whose non-trivial brackets are $[e_2,e_4]=e_1, [e_3,e_5]=e_1, [e_6,e_2]=e_4, [e_6,e_3]=e_5$. 
	
	For $t=\pi$, the rational basis $\B=\{f_1,\ldots,f_5\}=\{\pi e_1, -\pi e_4, e_2, -\pi e_5, e_3\}$ of $\n_{5,3}$ satisfies that the only non-trivial Lie brackets are $[f_2,f_3]=f_1$ and $[f_4,f_5]=f_1$. Moreover, we have that $[\phi(\pi)]_\B=(1)\oplus\smatriz{-1&1\\0&-1}^{\oplus 2}\in\SL(5,\Z)$. By Theorem \ref{theorem: yamada}, $\Gamma=\pi \Z\ltimes_\phi \exp^{N_{5,3}} \text{span}_\Z \B$ is a lattice in $S_{6,167}$. Thus, the corresponding complex solvmanifold has holomorphically trivial canonical bundle. 
	
	\medskip
	
	
	$\bullet$  $\s_{6,217}^{0,0,-1,-1}=(e^{16}-e^{25},e^{15}+e^{26},-e^{36},-e^{46},0,0)$ can be written as $(\R e_5 \times \R e_6)\ltimes_{\ad} \R^4$, where $\ad e_5|_{\R^4}=\smatriz{0&-1\\1&0}\oplus (0)^{\oplus 2}$ and $\ad e_6|_{\R^4}=\operatorname{diag}(1,1,-1,-1)$. Hence, the corresponding simply connected Lie group is $S_{6,217}^{0,0,-1,-1}:=\R^2\ltimes_\phi \R^4$, where
	\[  \phi(t,0)=\matriz{\cos t&-\sin t\\ \sin t & \cos t}\oplus \I_2, \quad \phi(0,s)=\operatorname{diag}(e^s,e^s,e^{-s},e^{-s}).\] The non-zero closed $(3,0)$-form $\tau$ from \eqref{eq: tau} is $\tau=\exp(it)(e^1+ie^2)\wedge (e^3+ie^4) \wedge (e^5+ie^6)$. Hence, we must choose $t=2\pi$. The corresponding lattice is a lattice in $S_{5,9}^{1,-1,-1}\times \R$, which admits a non-vanishing left-invariant closed $(3,0)$-form with respect to the same complex structure. As in previous cases, it can be seen that the corresponding solvmanifolds with holomorphically trivial canonical bundle are biholomorphic.  
	
	\medskip
	 
	 $\bullet$  $\s_{6,226}^{0,\beta,-1}=(e^{16}+e^{25},-e^{15}+e^{26},-e^{36}+\beta e^{45},-\beta e^{35}-e^{46},0,0)$, where $0<\beta< 1$, can be expressed as $(\R e_5 \times \R e_6)\ltimes_{\ad} \R^4$, where $\ad e_5|_{\R^4}=\smatriz{0&1\\-1&0}\oplus \smatriz{0&\beta\\-\beta&0}$ and $\ad e_6|_{\R^4}=\operatorname{diag}(1,1,-1,-1)$. The corresponding simply connected Lie group is $S_{6,226}^{0,\beta,-1}:=\R^2\ltimes_\phi \R^4$, where
	 \[ \phi(t,0)=\matriz{\cos t&\sin t\\ -\sin t&\cos t}\oplus \matriz{\cos(\beta t)&\sin(\beta t)\\ -\sin(\beta t)&\cos(\beta t)}, \quad \phi(0,s)=\operatorname{diag}(e^s,e^s,e^{-s},e^{-s}).\] The closed $(3,0)$-form $\tau$ from \eqref{eq: tau} is $\tau=\exp(-i(1-\beta)t)(e^1-ie^2)\wedge(e^3+ie^4)\wedge(e^5+ie^6)$.
	 As in previous cases we are forced to choose $t=\frac{2\pi k}{1-\beta}$ , for some $k \in \Z\setminus\{0\}$, ensuring that $\beta t\in t+2\pi \Z$. Therefore, the corresponding lattices (for the countable values of $\beta$ such that a lattice exists) can also be obtained considering the simply connected Lie group associated to $\s_{6,226}^{0,1,-1}$, which admits a non-vanishing left-invariant closed $(3,0)$-form respect to the same complex structure. Consequently, it can be checked that the corresponding complex solvmanifolds with holomorphically trivial canonical bundle are biholomorphic. 
	 
	 To get lattices we can choose, for instance, $\beta=\frac{1}{n}, n\in \N, n\geq 2$. Then, for $t=2\pi n$ and $s=t_m$, with $t_m$ defined as for $S_{5,11}^{\alpha,\alpha,-\alpha}\times \R$, it is easy to see that there is a rational basis satisfying the hypotheses of Theorem \ref{theorem: yamada}.
	 \medskip
	 	 
	$\bullet$ For the Lie algebra \[\s_{6,228}^{-\beta,\beta,\gamma,-\gamma}=(-\beta e^{15}+\gamma e^{16}+e^{25},-e^{15}-\beta e^{25}+\gamma e^{26},\beta e^{35}-\gamma e^{36}+e^{46}, -e^{36}+\beta e^{45}-\gamma e^{46},0,0),\] where $0\leq \beta\leq \gamma$, $0<\gamma$, we need to change the basis in order to apply Proposition \ref{proposition: tau explicit}. We take the ordered basis $\B=\{e_5+e_6,e_6-e_5,e_1,e_2,e_3,e_4\}$. Let $f_i$ denote the elements of $\B$, $1\leq i\leq 6$. Then, the complex structure defined by $Je_{2i-1}=e_{2i}$, $1\leq i\leq 3$ also satisfies $Jf_{2i-1}=f_{2i}$. Furthermore, we can express the Lie algebra as $(\R f_1\times \R f_2)\ltimes_{\ad} \R^4$, where \[\ad f_1|_{\R^4}=\smatriz{\gamma-\beta&1\\-1&\gamma-\beta}\oplus \smatriz{\beta-\gamma&1\\-1&\beta-\gamma} \quad \text{and} \quad \ad f_2|_{\R^4}=\smatriz{\gamma+\beta&-1\\1&\gamma+\beta}\oplus \smatriz{-\gamma-\beta&1\\-1&-\gamma-\beta}.\] The associated simply connected Lie group is $S_{6,228}^{-\beta,\beta,\gamma,-\gamma}:=\R^2\ltimes_\phi \R^4$, where 
	 \[ \phi(t,0)= e^{-t(\beta-\gamma)} \matriz{\cos t&\sin t\\-\sin t& \cos t}\oplus e^{t(\beta-\gamma)}\matriz{\cos t&\sin t\\-\sin t& \cos t}, \quad \text{and}\]\[ \phi(0,s)=e^{-s(\beta+\gamma)}\matriz{\cos s&-\sin s\\\sin s&\cos s}\oplus e^{s(\beta+\gamma)} \matriz{\cos s&\sin s\\-\sin s& \cos s}.\]
	 The closed $(3,0)$-form $\tau$ from \eqref{eq: tau} is  $\tau=\exp(-2it)(f^1+if^2)\wedge (f^3+if^4)\wedge (f^5+if^6)$. To construct a lattice $\Gamma$ using Theorem \ref{theorem: yamada} such that $\tau$ is $\Gamma$-invariant we must select $t\in\{2\pi, \pi\}$. For $t=2\pi$, choosing $\beta=\gamma$, for certain values of $s$ and $\gamma$ we can obtain lattices. However, these lattices can also be obtained considering the simply connected Lie group associated to the Lie algebra $\s_{5,13}^{2\gamma,-2\gamma,1}\times \R$, after changing to the basis $\{e_6,e_5,e_2,e_1,e_3,e_4\}$, which admits a non-vanishing left-invariant closed $(3,0)$-form with respect to the left-invariant complex structure we considered for $S_{6,228}^{-\beta,\beta,\gamma,-\gamma}$. The associated complex solvmanifolds with holomorphically trivial canonical bundle will be biholomorphic. 
	 
	For $t=\pi$, we set $\beta=\gamma=\frac{t_m}{2s}$, where $t_m=\log \left( \frac{m+\sqrt{m^2-4}}{2}\right)$, for $m\geq 3$. For $s=2\pi$, the (rational) basis $\{e_1+e_3, e^{t_m}e_1+e^{-t_m}e_3, e_2+e_4, e^{t_m}e_2+e^{-t_m}e_4\}$ of $\R^4$ satisfies that $[\phi(\pi,0)]_\B=-\I_4\in \SL(4,\Z)$ and $B_m:=[\phi(0,2\pi)]_\B=\smatriz{0&-1\\1&m}^{\oplus 2}\in\SL(4,\Z)$. According to Theorem \ref{theorem: yamada}, $\Gamma_m=(\pi \Z\oplus 2\pi\Z)\ltimes_\phi \text{span}_\Z \B$ is a lattice in $S_{6,228}^{-\beta,\beta,\gamma,-\gamma}$. Thus, the corresponding complex solvmanifold has holomorphically trivial canonical bundle.
	

\bigskip

In light of our computations, we have the following classification of the 6-dimensional solvable (non-nilpotent) Lie algebras corresponding to complex solvmanifolds with holomorphically trivial canonical bundle.

\begin{theorem}\label{theorem: classification}
	Let $(\Gamma\backslash G,J)$ be a $6$-dimensional complex solvmanifold with holomorphically trivial canonical bundle. Then the Lie algebra $\g$ of $G$ is isomorphic to $\g_i$, $1\leq i\leq 10$, or to one of the Lie algebras described in Theorem \ref{theorem: CRF}. 
\end{theorem}

As a consequence of the classification we can state the following result about completely solvable 6-dimensional Lie algebras. Recall that a solvable Lie algebra $\g$ is completely solvable if the adjoint operators $\ad x:\g\to\g$, with $x\in \g$, have only real eigenvalues.

\begin{theorem}\label{theorem: CS}
	Let $\g$ be a $6$-dimensional completely solvable Lie algebra equipped with a complex structure $J$ such that the corresponding simply connected $G$ admit a lattice $\Gamma$. If $K_{(\Gamma\backslash G,J)}$ is trivial, then there exists a non-zero holomorphic $(3,0)$-form defined on $\g$.
\end{theorem}
\begin{proof}
	According to Theorem \ref{theorem: classification}, since the complex solvmanifold has holomorphically trivial canonical bundle, we must have that the Lie algebra is isomorphic to $\g_i$, $1\leq i\leq 10$, or to one of the Lie algebras in Theorem \ref{theorem: CRF}. Since the only completely solvable Lie algebras among those are $\g_1$ and $\g_{5}$, the statement follows.
\end{proof}

To conclude, we provide next an example of a complex solvmanifold with holomorphically trivial canonical bundle which is not biholomorphic to one with an invariant holomorphic section of its canonical bundle.

\begin{proposition}\label{prop: S6,154}
	The complex solvmanifold with holomorphically trivial canonical bundle $(\Gamma\backslash S_{6,154}^0, J)$ previously constructed is not biholomorphic to a complex solvmanifold with an invariant holomorphic section of its canonical bundle.
\end{proposition}

\begin{proof}
	First, note that $M=\Gamma\backslash S_{6,154}^0$ is not even homeomorphic to a nilmanifold since the lattice $\Gamma$ is not nilpotent. This can be seen for instance from the identity $t_1^{t_6}=t_1^{-1}$. Moreover, we claim that $(M,J)$ is not biholomorphic to a complex solvmanifold arising from any of the Lie algebras $\g_1,\ldots,\g_{10}$. Indeed, given that $b_1(\Gamma\backslash S_{6,154}^0)=1$,  $b_1(\g_1)=b_1(\g_2^{\alpha})=b_1(\g_3)=b_1(\g_8)=b_1(\g_{10})=2$, together with the fact that $b_1(\g)\leq b_1(\Gamma\backslash G)$ for any solvmanifold $\Gamma\backslash G$, we deduce that $\Gamma\backslash S_{6,154}^0$ is not even homeomorphic to a solvmanifold arising from the algebras $\g_1, \g_2^{\alpha}, \g_3, \g_8$ or $\g_{10}$. Furthermore, according to 
	\cite[Theorem 5.1, Theorem 7.2, Theorem 8.1]{FPAN2}, the Lie algebra $\s_{6,154}^0$ does not admit a complex structure compatible with either an SKT, an LCK or a balanced metric, and from \cite[Theorem 1]{MM} (where the notation $\g_{6,83}^{0,0}$ is used), it follows that $\s_{6,154}^0$ does not admit a symplectic form. In contrast, according to \cite[Theorem 4.1, Remark 4.2, Theorem 4.5]{FOU}, any complex structure admitting a non-zero holomorphic $(3,0)$-form on $\g_4$ (resp. $\g_5, \g_7$) admits SKT (resp. balanced) metrics. Thus, if $M$ were biholomorphic to a complex solvmanifold $(\tilde{\Gamma}\backslash G_i, \tilde{J})$, for some $i=4,5,7$ then $M$ would admit a complex structure compatible with an SKT or a balanced metric. By applying Belgun's symmetrization process one would obtain a left-invariant  balanced or SKT metric on $S_{6,154}^0$ (\cite[Theorem 7]{Belgun}, \cite[Theorem 2.1]{FG}, see also \cite[Lemma 4.1]{P}), contradicting the fact that $\s_{6,154}^0$ does not admit any of these structures. For $\g_9\simeq \s_{6,152}$, given that $\g_9$ admits a symplectic form, as noticed in \cite[Theorem 1]{MM} (where it is referred to as $\n_{6,84}^1$), a diffeomorphism of $M$ with a solvmanifold arising from $\g_9$ would produce a symplectic form on $M$, which is again a contradiction due to Remark \ref{rmk: 154 no symplectic}.
	
	Finally, if $M$ were biholomorphic to a complex solvmanifold arising from $\g_6\simeq \s_{6,166}^{-1}$, then $(M,J)$ would admit a Vaisman metric, due to \cite[Proposition 3.6]{FOU} and \cite[Example 6.11]{AO}. However, since $b_1(\s_{6,154}^0)=b_1(\Gamma\backslash S_{6,154}^0)$, the argument in the proof of \cite[Theorem 4.3]{Ka} shows the existence of a left-invariant LCK metric on $S_{6,154}^0$, which contradicts the fact that $\s_{6,154}^0$ does not admit such a structure.
	
	This shows that $(M,J)$ cannot be biholomorphic to a complex solvmanifold with an invariant holomorphic section of its canonical bundle and finishes the proof. 
\end{proof}

\begin{remark}\label{rmk: candidatos}
	The complex solvmanifolds associated to $\s_{5,8}^0 \times \R$, $\s_{6,44}$, $\s_{6,145}^0$, $\s_{6,165}^{\alpha}$, $\s_{6,167}$ and $\s_{6,228}^{-\beta,\beta,\gamma,-\gamma}$, which we constructed using the value $\pi$, can be all seen to be non-homeomorphic to a nilmanifold, as the corresponding lattices are not nilpotent. Moreover, with similar arguments to those used for $\s_{6,154}^0$, we can rule out the biholomorphism with a complex solvmanifold associated to some of the Lie algebras $\g_i$, $1\leq i\leq 10$, though not all of them. For instance, since the complex solvmanifold $(\Gamma\backslash S_{6,44},J)$ can be seen to have first Betti number equal to $1$ and $\s_{6,44}$ does not admit a complex structure with either a balanced, an LCK or an SKT metric (\cite[Theorem 5.1, Theorem 7.2, Theorem 8.1]{FPAN2}), if the complex solvmanifold we constructed were biholomorphic to one with an invariant holomorphic section, it would have to be associated with the Lie algebra $\g_9$. It is worth pointing out that classifying all the lattices in a solvable Lie group up to isomorphism is a task that is generally not feasible, except in some special cases. Nevertheless, the complex solvmanifolds with holomorphically trivial canonical bundle that we have constructed offer interesting examples, as they are not readily identifiable as biholomorphic to a complex solvmanifold with an invariant holomorphic section of the canonical bundle. 
\end{remark}
	
	\newpage 
	
\section{Appendix: Tables} \label{S: algebras}
We present here the list of all $6$-dimensional solvable (non-nilpotent) strongly unimodular Lie algebras, together with their nilradical, using the naming convention of \cite{SW}. The column labeled \comillas{$J$} indicates whether the Lie algebra admits a complex structure or not. If there are specified parameters, it means that the Lie algebra admits a complex structure for those values, and does not admit for the other values. Finally, the column labeled \comillas{$d\psi=0$} indicates whether the Lie algebra admits a complex structure $J$ such that the associated Koszul $1$-form is closed. The symbol $\doublecheck$ indicates that there exists a complex structure $J$ with $\psi\equiv 0$.
\begin{table}[H]
	\renewcommand{\arraystretch}{1.1}
	\centering
\small{	\begin{tabular}{|l|l|l|c|c|} 
		\hline
		$\nil(\g)$ & $\g$ & Structure equations  
		 & 
		$J$ & $d\psi=0$\\ [0.1ex] 
		\hline
		$\R^5$ & $\s_{3,1}^{-1} \times \R^3$ & $(e^{13},-e^{23},0,0,0,0)$  
		 &$\times$& $\times$  \\  \hline

		$\R^2 \times \h_3$ & $\s_{3,1}^{-1} \times \h_3$ & $(e^{13},-e^{23},0,-e^{56},0,0)$ 
		& $\times$& $\times$\\  \hline

	    $\R^4$ & $\s_{3,1}^{-1} \times \s_{3,1}^{-1}$ & $(e^{13}, -e^{23}, 0, e^{46}, -e^{56}, 0)$ 
	    & $\times$ & $\times$ \\ \hline

	    $\R^5$ & $\s_{3,3}^0 \times \R^3$ & $(e^{23},-e^{13},0,0,0,0)$ 
	    & $\checkmark$ & $\checkmark$ \\ 
	    
	    \hline

	    $\R^2\times \h_3$ & $\s_{3,3}^0 \times \h_3$ & $(e^{23},-e^{13},0,-e^{56},0,0)$ 
	    & $\checkmark$ & $\checkmark$ \\ \hline

	    $\R^4$ & $\s_{3,3}^0 \times \s_{3,3}^0$ & $(e^{23},-e^{13},0,e^{56},-e^{46},0)$ 
	    & $\checkmark$ & $\checkmark$ \\ \hline

	    $\R^4$ & $\s_{3,1}^{-1} \times \s_{3,3}^0$ & $(e^{13},-e^{23},0,e^{56},-e^{46},0)$ 
	    & $\checkmark$ & $\doublecheck$
		\\
		\hline
	\end{tabular}}
	\caption{\textbf{Table 1.} Lie algebras $\g\times \h$, $\dim \g=\dim \h=3$.}
	\label{table:Decomposable 3x3}
\end{table}

	\begin{table}[H]
		\renewcommand{\arraystretch}{1.2}
	\small{	\centering
		\begin{tabular}{|l|l|l|c|c|} 
			\hline
			$\nil(\g)$ & $\g$ & Structure equations  
			 & $J$ & $d\psi=0$\\ 
			\hline
			
			\multirow{3}*{$\R^5$} & $\s_{4,3}^{a,-(1+a)}\times \R^2$ {\tiny($-1<a\leq -\frac12$)}  & \multirow{1}*{$(e^{14},ae^{24},-(1+a)e^{34},0,0,0)$}  
			 & \multirow{1}*{{\tiny$a=-\frac12$}} & \multirow{1}*{$\times$}
		\\ \cline{2-5}

			&$\s_{4,4}^{-2}\times \R^2$ & $(e^{14}+e^{24},e^{24},-2e^{34},0,0,0)$ 
			& $\times$ & $\times$ \\ \cline{2-5}&&&& \\ [-1.2em]

			&$\s_{4,5}^{\alpha,-\frac12 \alpha}\times \R^2$ {\tiny($\alpha>0$)} & \multirow{1}*{$(\alpha e^{14},-\frac{\alpha}{2} e^{24}+e^{34},-e^{24}-\frac{\alpha}{2}e^{34},0,0,0)$} 
			&\multirow{1}*{$\checkmark$} & \multirow{1}*{$\times$} \\ \hline

			$\h_3 \times \R^2$ & $\s_{4,6}\times \R^2$ & $(-e^{23}, e^{24}, -e^{34},0,0,0)$ 
			& $\checkmark$ & $\times$ \\ \cline{2-5}
			
			& $\s_{4,7}\times\R^2$ &  $(-e^{23}, e^{34}, -e^{24},0,0,0)$ 
			& $\checkmark$ & $\checkmark$ 
			\\
			\hline
		\end{tabular}}
		\caption{\textbf{Table 2.} Lie algebras $\g\times \R^2$, $\dim \g=4$, $\g$ indecomposable.}
		\label{table:Decomposable 4+R2}
	\end{table}
	
	\begin{table}[H]
	\renewcommand{\arraystretch}{1.1}
	\centering
{\small		\begin{tabular}{|l|l|c|c|} 
			\hline
			$\g$  & Structure equations &
			$J$ & $d\psi=0$\\  
			\hline 
			
			$\s_{5,3}^{-1}\times \R$  & $(e^{25},0, e^{35}, -e^{45}, 0,0)$ & $\times$ & $\times$ \\ \hline 
			
			$\s_{5,4}^0\times \R$  & $(e^{25},0,e^{45},-e^{35},0,0)$
			& $\checkmark$ & $\checkmark$ \\ \hline

			$\s_{5,6}^{-1}\times \R$ & $(e^{15}+e^{25},e^{25},-e^{35}+e^{45},-e^{45},0,0)$
			& $\times$ & $\times$ \\   \hline

			$\s_{5,7}^{-3}\times \R$  & $(e^{15}+e^{25},e^{25}+e^{35},e^{35},-3 e^{45},0,0)$ &  $\times$ & $\times$ \\ \hline

			$\s_{5,8}^0\times \R$  & $(e^{25}+e^{35},-e^{15}+e^{45},e^{45},-e^{35},0,0)$
			& $\checkmark$ & $\checkmark$ \\ \hline

			$\s_{5,9}^{a,b,-a-b-1}\times \R$   & \multirow{2}*{$(e^{15}, a e^{25}, b e^{35}, -(a+b+1)e^{45},0,0)$}& \multirow{2}*{{\tiny$(a,b)=(1,-1)$}} & \multirow{2}*{$\doublecheck$} \\ [-0.3em]
			{\tiny $0<|a+b+1|\leq |b|\leq |a|\leq 1$}&&&\\ \hline

			$\s_{5,10}^{a,-(a+2)}\times \R$  & \multirow{2}*{$(e^{15}+e^{25}, e^{25}, a e^{35}, -(a+2)e^{45},0,0)$}
			& \multirow{2}*{$\times$} & \multirow{2}*{$\times$} \\  [-0.3em]
			{\tiny$a\leq -1, a\neq -2$}
			&&&\\ \hline &&&\\ [-1.2em]

			$\s_{5,11}^{\alpha,\beta,-\frac{\alpha+\beta}{2}}\times \R$  & \multirow{2}*{$(\alpha e^{15}, \beta e^{25}, -\frac{\alpha+\beta}{2} e^{35}+e^{45}, -e^{35}-\frac{\alpha+\beta}{2} e^{45},0,0)$}
			& \multirow{2}*{{\tiny$\beta=\alpha$}} & \multirow{2}*{$\checkmark$} \\ [-0.3em]
			{\tiny $\alpha>0$, $\beta\neq0$, $|\beta|\leq\alpha$} &&&\\ \hline

			$\s_{5,12}^{-1,\beta}\times \R $ {\tiny($\beta>0$)}   & \multirow{1}*{$(e^{15}+e^{25},e^{25}, -e^{35}+\beta e^{45}, -\beta e^{35}-e^{45},0,0)$}&  \multirow{1}*{$\times$} & \multirow{1}*{$\times$} 
		\\ \hline 
			
			$\s_{5,13}^{\alpha,-\alpha,\gamma}\times \R$  &  $(\alpha e^{15}+e^{25}, -e^{15}+\alpha e^{25}, $
			& \multirow{2}*{$\checkmark$} & $\gamma\neq1\;\checkmark$    \\ [-0.3em]
			{\tiny $0<\gamma\leq 1$, $0\leq \alpha$} &\quad$-\alpha e^{35}+\gamma e^{45}, -\gamma e^{35}-\alpha e^{45}, 0 ,0)$&& $\gamma=1\;\doublecheck$\\
			\hline \end{tabular}
		\caption{\textbf{Table 3.} Lie algebras $\g\times \R$, $\dim \g=5$, $\g$ indecomposable and $\nil(\g)=\R^5$.
		}
		\label{table:Almost abelian 5+R}}
	\end{table} 
			 
		\begin{table}[H]
		 	\renewcommand{\arraystretch}{1.2}
		 	\centering
		 		\begin{tabular}{|l|l|l|c|c|} 
		 		\hline
		 		$\nil(\g)$ & $\g$ & Structure equations & $J$ & $d\psi=0$\\ [0.5ex] 
			 		\hline

			 		\multirow{2}*{$\h_3 \times \R^2$} & $\s_{5,15}\times\R$ & $(-e^{23}+e^{45}, e^{25}, -e^{35}, 0,0,0)$
			 		& $\times$ & $\times$\\ \cline{2-5}

			 		& $\s_{5,16}\times \R$ &  $(-e^{23}+e^{45}, e^{35}, -e^{25}, 0 ,0 ,0)$	& $\checkmark$ & $\checkmark$ \\ \hline

			 \multirow{2}*{$\R^4$} & $\s_{5,41}^{-1,-1}\times \R$ & $(e^{14}, e^{25}, -e^{34}-e^{35}, 0 , 0, 0)$
			 & $\times$ & $\times$ \\ \cline{2-5}
			 
			 & $\s_{5,43}^{-2,0}\times \R$ & $(-2e^{14}, e^{24}+e^{35}, -e^{25}+e^{34},0,0,0)$
			 & $\checkmark$ & $\times$ \\
			\hline
		\end{tabular}
				\caption{\textbf{Table 4.} Lie algebras $\g\times \R$, $\dim \g=5$, $\g$ indecomposable and $\nil(\g)\neq \R^5$.}
		\label{table:Decomposable 5+R}
	\end{table}

\begin{table}[h]
	\renewcommand{\arraystretch}{1.3}
	\centering
	{\small\begin{tabular}{|l|l|c|c|} 
		\hline
		$\g$ & Structure equations & $J$ & $d\psi=0$\\ [0.5ex] 
		\hline
		
		$\s_{6,4}^{-1}$ & $(e^{26},e^{36},0,e^{46},-e^{56},0)$
		&$\times$ & $\times$ \\ \hline

		$\s_{6,5}^0$ & $(e^{26},e^{36},0,e^{56},-e^{46},0)$
		&$\times$ & $\times$ \\ \hline&&&\\ [-1.5em]

		$\s_{6,7}^{-\frac12}$  & $(-\tfrac12 e^{16}+e^{26}, -\tfrac12 e^{26}, e^{46}, 0, e^{56}, 0)$ & $\times$ & $\times$ \\ \hline

		$\s_{6,8}^{a,-(a+1)}$  & \multirow{2}*{$(e^{26},0,e^{36},a e^{46}, -(a+1)e^{56},0)$} & \multirow{2}*{$\times$} & \multirow{2}*{$\times$} \\
		$-1<a\leq -\tfrac12$ &&& \\ \hline 
		
		$\s_{6,9}^{\alpha,-\frac{\alpha}{2}}$& \multirow{2}*{$(e^{26},0,\alpha e^{36}, -\frac{\alpha}{2}e^{46}+e^{56}, -e^{46}-\frac{\alpha}{2}e^{56},0)$} & \multirow{2}*{$\times$} & \multirow{2}*{$\times$} \\
		$\alpha>0$ &&& \\ \hline&&&\\ [-1.5em]

		$\s_{6,11}^{-\frac{3}{2}}$  & $(e^{16}+e^{26},e^{26}+e^{36},e^{36},-\frac32 e^{46}+e^{56}, -\frac32 e^{56},0)$ & $\times$ & $\times$  \\ \hline&&&\\ [-1.5em] 
		
		\multirow{2}*{$\s_{6,12}^{-\frac14}$}  & $(-\frac14 e^{16}+e^{26}, -\frac14 e^{26}+e^{36}, -\frac14 e^{36}+e^{46},$& \multirow{2}*{$\times$} & \multirow{2}*{$\times$} \\ 
		&$-\frac14 e^{46}, e^{56},0)$ && \\ \hline

		$\s_{6,13}^{a,-(3a+1)}$  & \multirow{2}*{$(a e^{16}+e^{26}, a e^{26}+e^{36}, a e^{36}, e^{46}, -(3a+1) e^{56},0)$}  & \multirow{2}*{$\times$} & \multirow{2}*{$\times$} \\
		$a\in [-\tfrac23,0)\setminus\{-\tfrac13\}$ &&&\\ \hline&&&\\ [-1.5em]

		$\s_{6,14}^{a,-(a+\frac12)}$ & $(a e^{16}+e^{26}, a e^{26}, $-{\tiny$(a+\frac12)$}$e^{36}+e^{46},$   & \multirow{2}*{{\tiny$a=-\frac14$}} & \multirow{2}*{$\times$} \\
		$a\leq -\frac14, a\neq -\frac12$ & $-${\tiny$(a+\frac12)$}$e^{46}, e^{56}, 0)$ &&\\ \hline&&&\\ [-1.5em]

		$\s_{6,15}^{\alpha, -\frac32 \alpha}$ & $(\alpha e^{16}+e^{26}, \alpha e^{26}+e^{36}, \alpha e^{36},$ & \multirow{2}*{$\times$} & \multirow{2}*{$\times$} \\
		$\alpha>0$ &\quad$-\frac{3\alpha}{2}e^{46}+e^{56}, -e^{46}-\frac{3\alpha}{2} e^{56},0)$&&\\ \hline 	
		
		$\s_{6,16}^{\alpha,-4\alpha}$ & $(\alpha e^{16}+e^{26}+e^{36}, -e^{16}+\alpha e^{26}+e^{46}$, & \multirow{2}*{$\checkmark$} & \multirow{2}*{$\times$} \\
		$\alpha<0$ & $\quad\alpha e^{36}+e^{46}, -e^{36}+\alpha e^{46}, -4\alpha e^{56},0)$  &&\\ \hline

		$\s_{6,17}^{a,b,c,d}$, {\tiny$a+b+c+d=1$} & \multirow{2}*{$(e^{16}, a e^{26}, b e^{36}, c e^{46}, d e^{56},0)$} & \multirow{2}*{{\tiny$a=1, c=b$}} & \multirow{2}*{$\times$} \\
		{\tiny $0<|d|\leq|c|\leq|b|\leq|a|\leq 1$} &&& \\ \hline

		$\s_{6,18}^{a,b,-2a-b-1}$  & \multirow{2}*{$(a e^{16}+e^{26}, a e^{26}, e^{36}, b e^{46}, -(2a+b+1)e^{56},0)$}& \multirow{2}*{{\tiny$(a,b)=(1,-\frac{3}{2})$}} & \multirow{2}*{$\times$} \\
		{\tiny $0<|2a+b+1|\leq|b|\leq 1, a\neq 0$} &&&\\ \hline&&&\\ [-1.5em]

		$\s_{6,19}^{\alpha,\beta,\gamma,-\frac{\alpha+\beta+\gamma}{2}}$ & $(\alpha e^{16}, \beta e^{26}, \gamma e^{36}, -${\tiny$\frac{\alpha+\beta+\gamma}{2}$}$ e^{46}+e^{56},$ & \multirow{2}*{{\tiny$\beta=\alpha$}} & \multirow{2}*{$\times$} \\
		$0<|\gamma|\leq|\beta|\leq\alpha$ &  $-e^{46}-${\tiny$\frac{\alpha+\beta+\gamma}{2}$}$ e^{56},0)$&&\\ \hline&&&\\ [-1.5em]
		
		$\s_{6,20}^{\alpha,\beta,-\alpha-\frac{\beta}{2}}$ & $(\alpha e^{16}+e^{26}, \alpha e^{26}, \beta e^{36}$,  & \multirow{2}*{{\tiny$\beta=\alpha$}} & \multirow{2}*{$\times$} \\
		$0<\alpha, 0\neq \beta$ &\quad$-(\alpha+\frac{\beta}{2})e^{46}+e^{56}, -e^{46}-(\alpha+\frac{\beta}{2})e^{56},0)$&&\\ \hline

		$\s_{6,21}^{\alpha,\beta,\gamma,-2\alpha-2\beta}, \;${\tiny $\alpha<-\beta,$}  & $(\alpha e^{16}+e^{26}, -e^{16}+\alpha e^{26}, \beta e^{36}+\gamma e^{46},$ & \multirow{2}*{$\checkmark$} & \multirow{2}*{$\times$} \\
		{\tiny $0<\gamma\leq 1$. If $\gamma=1$ then $\alpha\leq \beta$} &$\quad -\gamma e^{36}+\beta e^{46}, -2(\alpha+\beta)e^{56},0)$&&\\
		\hline
	\end{tabular}}
	\caption{\textbf{Table 5.} Indecomposable Lie algebras with $\nil(\g)=\R^5$.}
	\label{table:Almost abelian indecomposable}
\end{table}

	\begin{table}[H]
		\renewcommand{\arraystretch}{1.2}
		\centering
		\begin{tabular}{|l|l|l|c|c|} 
				\hline
				$\nil(\g)$ & $\g$ & Structure equations 
				& $J$ & $d\psi=0$\\  
				\hline
				
				\multirow{3}*{$\h_3\times \R^2$} & $\s_{6,24}$  & $(-e^{23},e^{26},-e^{36},e^{56},0,0)$ & $\times$ & $\times$ \\ \cline{2-5}

				& $\s_{6,25}$  & $(-e^{23},e^{36},-e^{26},0,e^{46},0)$  & $\checkmark$ & $\checkmark$ \\ \cline{2-5}
				
				& $\s_{6,30}$  & $(-e^{23}+e^{56},e^{26},-e^{36},0,e^{46},0)$  & $\times$ & $\times$ \\ \cline{2-5}

				& $\s_{6,31}$  & $(-e^{23}+e^{56},e^{36},-e^{26},0,e^{46},0)$ & $\times$ & $\times$ \\ \cline{2-5}

				& $\s_{6,32}^{-1}$  & $(-e^{23},e^{36},0,e^{46},-e^{56},0)$&  $\times$ & $\times$ \\ \cline{2-5}

				& $\s_{6,34}^0$ &  $(-e^{23},e^{36},0,e^{56},-e^{46},0)$ & $\times$ & $\times$ \\ \cline{2-5} 
				
				& $\s_{6,43}$  & $(-e^{23},e^{26},-e^{36},e^{26}+e^{46},e^{36}-e^{56},0)$ & $\times$ & $\times$ \\ \cline{2-5}

				& $\s_{6,44}$  & $(-e^{23},e^{36},-e^{26},e^{26}+e^{56},e^{36}-e^{46},0)$ & $\checkmark$ & $\checkmark$ \\ \cline{2-5} 
				
				& $\s_{6,45}^{a,-1}$ {\tiny$(a\neq0)$} & $(-e^{23}, a e^{26}, -a e^{36}, e^{46}, -e^{56},0)$  & $\times$ & $\times$ \\ 
				 \cline{2-5}

				& $\s_{6,46}^{\alpha,-\alpha}$ {\tiny($\alpha>0$)}  & $(-e^{23},e^{36},-e^{26},\alpha e^{46},-\alpha e^{56},0)$ & $\times$ & $\times$ \\ \cline{2-5}
				
				& $\s_{6,47}^{-1}$ & $(-e^{23},-e^{26},e^{36},e^{36}+e^{46},-e^{56},0)$ & $\times$ & $\times$ \\ \cline{2-5}
				
				& $\s_{6,51}^{\alpha,0}$ {\tiny ($\alpha>0$)} & $(-e^{23},\alpha e^{26},-\alpha e^{36},e^{56},-e^{46},0)$ & $\checkmark$ & $\times$ \\ \cline{2-5}
				
				& \multirow{2}*{$\s_{6,52}^{0,\beta}$ {\tiny($\beta>0$)}} & \multirow{2}*{$(-e^{23},e^{36},-e^{26},\beta e^{56},-\beta e^{46},0)$} & \multirow{2}*{$\checkmark$} & $\beta\neq 1\, \checkmark$ \\
				& & & & $ \beta=1\, \doublecheck$\\ \hline

				$\n_{5,1}$ & $\s_{6,140}^{-1}$ &$(e^{16}-e^{35},-e^{26}-e^{45},e^{36},-e^{46},0,0)$ & $\times$ & $\times$\\ \cline{2-5}

				& $\s_{6,145}^0$  & $(e^{26}-e^{35},-e^{16}-e^{45},e^{46},-e^{36},0,0)$& $\checkmark$ & $\checkmark$ \\ \cline{2-5}

				& $\s_{6,146}^{-1}$ & $(e^{16}-e^{35}+e^{36},-e^{26}-e^{45}-e^{46},e^{36},-e^{46},0,0)$& $\times$ & $\times$ \\ \cline{2-5}

				& $\s_{6,147}^0$ & $(e^{26}-e^{35},-e^{16}+e^{36}-e^{45},e^{46},-e^{36},0,0)$ & $\checkmark$ & $\doublecheck$ (New) \\ \hline

				$\n_{5,2}$ & $\s_{6,151}$ & $(e^{16}-e^{35},-e^{26}-e^{34}+e^{46},-e^{45},-e^{46},e^{56},0)$ & $\times$ & $\times$ \\ \cline{2-5}

				& $\s_{6,152}$ & $(-e^{26}-e^{35},e^{16}-e^{34}+e^{56},-e^{45},e^{56},-e^{46},0)$ & $\checkmark$ & $\doublecheck$ \\ \cline{2-5}

				& $\s_{6,154}^0$ & $(-e^{26}-e^{35},e^{16}-e^{34},-e^{45},e^{56},-e^{46},0)$&$\checkmark$ & $\checkmark$ \\ \cline{2-5}

				& $\s_{6,155}^{-1}$ & $(-e^{16}-e^{35},e^{26}-e^{34},-e^{45},e^{46},-e^{56},0)$& $\times$ &$\times$ \\ \hline

				$\n_{5,3}$ & $\s_{6,158}$ & $(-e^{24}-e^{35},0,e^{36},0,-e^{56},0)$ & $\checkmark$ & $\times$ \\ \cline{2-5}

				& $\s_{6,159}$ & $(-e^{24}-e^{35},0,-e^{56},0,e^{36},0)$ & $\checkmark$ & $\checkmark$ \\ \cline{2-5}

				& $\s_{6,160}$ & $(-e^{24}-e^{35},e^{46},e^{36},0,-e^{56},0)$ & $\times$ & $\times$ \\ \cline{2-5}

				& $\s_{6,161}^{\pm 1}$ & $(-e^{24}-e^{35}, \pm e^{46}, -e^{56},0,e^{36},0)$ & $\times$ & $\times$ \\ \cline{2-5}

				& $\s_{6,162}^a$ & \multirow{2}*{$(-e^{24}-e^{35},e^{26},a e^{36}, -e^{46}, -a e^{56},0)$} & \multirow{2}*{$a=1$} & \multirow{2}*{$\doublecheck$}\\ 
				& {\tiny$0<|a|\leq 1$} & & & \\ \cline{2-5}

				& $\s_{6,163}$ & $(-e^{24}-e^{35},e^{26},e^{26}+e^{36},-e^{46}-e^{56},-e^{56},0)$  & $\times$ & $\times$ \\ \cline{2-5}

				& $\s_{6,164}^\alpha$ {\tiny ($\alpha>0$)} & $(-e^{24}-e^{35},\alpha e^{26}, -e^{56}, -\alpha e^{46}, e^{36},0)$ & $\checkmark$ &$\times$\\ \cline{2-5}

				& \multirow{2}*{$\s_{6,165}^\alpha$ {\tiny ($\alpha>0$)}} & $(-e^{24}-e^{35},\alpha e^{26}-e^{36}, e^{26}+\alpha e^{36},$& \multirow{2}*{$\checkmark$} & \multirow{2}*{$\checkmark$}
				\\
				& & $\quad-\alpha e^{46}-e^{56}, e^{46}-\alpha e^{56},0)$&&\\ \cline{2-5}

				& $\s_{6,166}^{\alpha}$ & \multirow{2}*{$(-e^{24}-e^{35},-e^{46},-\alpha e^{56}, e^{26}, \alpha e^{36},0)$} & \multirow{2}*{$\checkmark$} & $\alpha\neq \pm 1\, \checkmark$  \\
				& {\tiny$0<|\alpha|\leq 1$} & & & $\alpha=\pm 1\, \ \doublecheck$\\ \cline{2-5}

				& $\s_{6,167}$ & $(-e^{24}-e^{35},-e^{36},e^{26},e^{26}-e^{56},e^{36}+e^{46},0)$ & $\checkmark$ & $\checkmark$ \\
				\hline
			\end{tabular}
		\caption{\textbf{Table 6.} Indecomposable Lie algebras with $\dim \nil(\g)=5$, $\nil(\g)\neq \R^5$.}
	\label{table:almost nilpotent}
\end{table}

 \begin{table}[H]
 	\renewcommand{\arraystretch}{1.2}
 	\centering
 {\small	\begin{tabular}{|l|l|c|c|} 
 			\hline
 			$\g$ & Structure equations & $J$ & $d\psi=0$\\ [0.5ex] 
 			\hline
 	 $\s_{6,204}^{-1,-1}$ & $(e^{56},-e^{25}-e^{26},e^{36},e^{45},0,0)$ & $\times$ & $\times$ \\ \hline

 	 $\s_{6,208}^{0,-2}$& $(e^{56},e^{26}-e^{35},e^{25}+e^{36},-2 e^{46},0,0)$& $\times$ & $\times$ \\ \hline

 	 $\s_{6,213}^{a,-a-1,c,-c-1}$ & $(ae^{15}+ce^{16},-(a+1)e^{25}-(c+1)e^{26},$ & {\tiny$(a,c) =(-\tfrac12,-\tfrac12),$} & \multirow{2}*{$\times$} \\ 
 	 $a\leq-\frac12, \quad c\in\R$ &$e^{36},e^{45},0,0)$&{\tiny$(-1,-2),(-2,-1)$}&\\ \hline

 	 $\s_{6,214}^{a,-(a+2),-1}$ & $(e^{15}+e^{26},e^{25},a e^{35}+e^{36},$ & \multirow{2}*{$\times$} & \multirow{2}*{$\times$} \\
 	 $a\in \R$ &$-(a+2)e^{45}-e^{46},0,0)$&&\\  \hline
 	 
 	$\s_{6,215}^{-2,-1}$ & $(e^{15}+e^{25},e^{25},e^{36},-2e^{45}-e^{46},0,0)$ & $\times$ & $\times$ \\ \hline

 	 $\s_{6,216}^{\alpha,-2\alpha,-1}$ & $(\alpha e^{15}+e^{25},-e^{15}+\alpha e^{25},e^{36},$& \multirow{2}*{$\checkmark$} & \multirow{2}*{$\times$} \\
 	 $\alpha<0$ &$-2\alpha e^{45}-e^{46},0,0)$&&\\ \hline

 	 $\s_{6,217}^{\alpha,-\alpha,\gamma,-(\gamma+2)}, $ {\tiny $\alpha>0$ or} & $(e^{16}-e^{25},e^{15}+e^{26},$& $\alpha> 0$ & $\times$ \\
 \cline{3-4}
 	 {\tiny $\alpha=0$ and $\gamma\geq-1, \gamma\neq0 $} &\quad$\alpha e^{35}+\gamma e^{36}, -\alpha e^{45}-(\gamma+2)e^{46},0,0)$& {\small$(\alpha,\gamma)=(0,-1)$} & $\checkmark$ \\ \hline

 	 $\s_{6,224}^{\alpha,-\alpha,-1}$ & $(\alpha e^{15}-e^{16}+e^{25},-e^{15}+\alpha e^{25}-e^{26},$ & \multirow{2}*{$\alpha > 0$} & \multirow{2}*{$\times$} \\ 
 	 $\alpha\geq 0$& \quad$-\alpha e^{35}+e^{36}+e^{46},-\alpha e^{45}+e^{46},0,0)$&&\\ \hline

 	 $\s_{6,226}^{0,\beta,-1}$ & $(e^{16}+e^{25},-e^{15}+e^{26},-e^{36}+\beta e^{45},$& \multirow{2}*{$\checkmark$}& $\beta\neq 1\; \checkmark$ \\ 
 	 $0<\beta\leq 1$&$-\beta e^{35}-e^{46},0,0)$&& $\beta=1\;\doublecheck
 	 $ \\ \hline

 	 $\s_{6,227}^{0,-1}$ & $(e^{16}+e^{25},-e^{15}+e^{26},-e^{36}+e^{45},-e^{46},0,0)$& $\checkmark$ & $\times$\\ \hline
 	 
$\s_{6,228}^{-\beta,\beta,\gamma,-\gamma}$ & $(-\beta e^{15}+\gamma e^{16}+e^{25},-e^{15}-\beta e^{25}+\gamma e^{26}$, & \multirow{2}*{$\checkmark$} & \multirow{2}*{$\checkmark$}\\
$0\leq \beta\leq \gamma, 0<\gamma$ 	 &\quad$\beta e^{35}-\gamma e^{36}+e^{46},-e^{36}+\beta e^{45}-\gamma e^{46},0,0)$ &  &\\ 
 	 \hline
  \end{tabular}}
\caption{\textbf{Table 7.} Indecomposable Lie algebras with $\nil(\g)=\R^4$.}
\label{table:nilradical codimension 2}
 \end{table}

\

\end{document}